\documentclass{article}

\usepackage{color}
\usepackage{array}

\usepackage[colorlinks,allcolors=blue]{hyperref}   

\usepackage[normalem]{ulem} 
\usepackage[margin=1in]{geometry}
\usepackage{url}   

\usepackage[ansinew]{inputenc}
\usepackage{fancyhdr}
\usepackage[active]{srcltx}
\usepackage{amsmath,amssymb}
\usepackage[pdftex]{graphicx}
\usepackage{authblk}
\usepackage{placeins}

\newcommand{\tmagenta}[1]{\textcolor{magenta}{#1}}

\newcommand{\AC}{\text{AC}}
\newcommand{\DC}{\text{DC}}

\newcommand{\vDL}{v^{\text{DL}}}
\newcommand{\uP}{u^{\text{P}}}
\newcommand{\uDL}{u^{\text{DL}}}

\newcommand{\shddot}{\!:\!}

\newcommand{\Rbb}{\mathbb{R}}

\newcommand{\cP}{{\cal P}}

\newcommand{\cS}{{\cal S}}

\newcommand{\micr}{\,$\mu$m\,}

\newcommand{\bfm}[1]{\mathbf{#1}}
\newcommand{\bsn}[1]{\boldsymbol{#1}}

\newcommand{\bfa}{\bfm{a}}

\newcommand{\bfb}{\bfm{b}}
\newcommand{\bfc}{\bfm{c}}

\newcommand{\bfu}{\bfm{u}}

\newcommand{\bfV}{\bfm{V}}

\newcommand{\bff}{\bfm{f}}

\newcommand{\bfF}{\bfm{F}}
\newcommand{\bfG}{\bfm{G}}
\newcommand{\bfU}{\bfm{U}}

\newcommand{\bfM}{\bfm{M}}

\newcommand{\bfK}{\bfm{K}}

\newcommand{\bfz}{\bfm{z}}

\newcommand{\bfC}{\bfm{C}}
\newcommand{\bfH}{\bfm{H}}

\newcommand{\bfphi}{\bsn{\phi}}
\newcommand{\bfPhi}{\bsn{\Phi}}

\newcommand{\bfPsi}{\bsn{\Psi}}
\newcommand{\bfUps}{\bsn{\Upsilon}}

\newcommand{\bfalpha}{\bsn{\alpha}}
\newcommand{\bfbeta}{\bsn{\beta}}
\newcommand{\bfgamma}{\bsn{\gamma}}

\newcommand{\bfmu}{\bsn{\mu}}

\begin{document}
\title{Virtual twins of nonlinear vibrating multiphysics microstructures: physics-based versus deep learning-based approaches}
\author[1]{Giorgio Gobat\footnote{giorgio.gobat@polimi.it}}
\author[2]{Stefania Fresca}
\author[2]{Andrea Manzoni}
\author[1]{Attilio Frangi}
\affil[1]{Department of Civil and Environmental Engineering, Politecnico di Milano, P.za Leonardo da Vinci 32, 20133 Milano, Italy}
\affil[2]{MOX - Department of Mathematics, Politecnico di Milano, P.za Leonardo da Vinci 32, 20133 Milano, Italy}

\maketitle
\section*{Abstract} 
Micro-Electro-Mechanical-Systems are complex structures, often involving nonlinearites of geometric and multiphysics nature, that are used as sensors and actuators in countless applications. Starting from full-order representations, we apply deep learning techniques to generate accurate, efficient and real-time reduced order models to be used as virtual twin for the simulation and optimization of higher-level complex systems.
We extensively test the reliability of the proposed procedures on micromirrors, arches and gyroscopes, also displaying intricate dynamical evolutions like internal resonances. In particular, we discuss the accuracy of the deep learning technique and its ability to replicate and converge to the invariant manifolds predicted using the recently developed direct parametrization approach that allows extracting the nonlinear normal modes of large finite element models.
Finally, by addressing an electromechanical gyroscope, we show that the non-intrusive deep learning approach generalizes easily to complex multiphysics problems.

\section{Introduction}
\label{sec:intro}

Even if the fast development of computational resources has enabled the simulation of complex structures, involving multi-scale and multi-physics phenomena, the ever growing computational costs continuously drives the search of efficient, but accurate, reduced order modeling techniques. Well-known strategies are nowadays routinely applied in linear problems arising from vibratory systems and make use of linear normal modes as an optimal projection basis upon which the equations of motion can be reduced. However, their extensions to problems in nonlinear dynamics is  still an open issue, and thus object of intensive research.
Early ideas to project the nonlinear equations of motion
onto a selected subset of the linear modes basis are progressively being abandoned. 
Indeed, the presence of strong nonlinear coupling terms between low frequency {\it master} modes and high-frequency ones makes the aforementioned approach not efficient, and the initial choice of the trial space is a delicate matter \cite{cyrilstep}.
The Proper Orthogonal Decomposition (POD) method offers a potential gain since it optimally orients the subspace bases to better fit the curvatures of the nonlinear data \cite{kerschen2005,amabili07,amabili2003reduced}, and enables the construction of such a subspace from data collected by simulating the behavior of the system for a restricted amount of configurations; however, the use of linear subspaces -- possibly of local type -- might still represent a computational bottleneck  
\cite{POD21}.
Other approaches, like the implicit condensation and expansion \cite{ijnm19} or modal derivatives methods, have been introduced with the aim of taking the amplitude dependence of modes
into account. However, they
are limited to moderate transformations 
and apply only in presence of sufficient slow/fast separation between the slave and master coordinates \cite{comparison21}. 

On the other hand, truly nonlinear methods for nonlinear vibratory systems are, so far, strictly connected to the concepts of Nonlinear Normal Modes (NNMs) and invariant manifolds, as keys to compute accurate Reduced Order Models (ROMs) \cite{SP91,SP93}. These start by defining a nonlinear relationship between the
original coordinates and those of the reduced dynamics.
In the conservative case, the existence and uniqueness of the searched invariant structures are framed by the Lyapunov centre theorem which guarantees, under non-resonance conditions, the existence of a smooth manifold, tangent at the origin to the associated linear eigenspace, which is invariant. 
This means that a trajectory of an autonomous system initiated on the manifold will develop on the manifold itself.
For dissipative systems, the picture gets more complicated, as the whole phase space is foliated by invariant manifolds tangent at the origin to the linear subspace \cite{haller18}. 
The application of these methods to large Finite Element Models (FEMs) has remained sporadic until recently, but is currently receiving an impressive boost by contributions \cite{vizza2020,NNM21,vizzaccaro2021high,jain2022compute}
in which direct approaches, called Direct Parametrization of Invariant Manifolds (DPIM), bypass the requirement of computing the whole modal basis. 
Applications to complex structures with millions of degrees of freedom (DOFs) and featuring also internal resonances and parametric excitation have been recently demonstrated in 
\cite{vizzaccaro2021high,opreni22}.
However, their extension to coupled problems and nonlinearities of generic type is still an open issue, and requires dedicated developments. 
In parallel, exploiting machine learning methods for constructing {\it real-time} surrogates, namely virtual twins or \tmagenta{surrogates}
of nonlinear dynamical systems, has become an
area of increasing interest for the system dynamics community. A virtual twin can be considered as a particularly accurate and fast reduced order model
of the physical asset, that might be used in the simulation and/or optimization of complex systems in which the asset is allocated or as a key component enabling digital twins.

Even if there is still no consolidated view on the definition of digital twins, they are often interpreted \cite{ritto2022transfer,ritto2021digital,booyse2020deep,agathos2021crack,bigoni2020simulation} as consisting of a physical asset, a virtual representation of that
product, and the data connections that feed data
from the physical to the virtual representation and vice-versa through 
mirroring or twinning \cite{jones,cimino,ganguli}.
While the operational feedback is naturally provided by the proposed virtual twin -- 
through the ability of computing quantities and outputs of interest -- 
the twinning stage, realized by means of a physical-to-virtual connection, is not considered in this contribution; nevertheless, it could be profitably included in the process. 

The idea of leveraging the modeling capability and learning flexibility of Deep Learning (DL) methods to
identify at once nonlinear transformations, nonlinear invariant manifolds and modal dynamics of NNMs from the system response data only, is appealing.
Among others, great success has been encountered by the Physics Informed Neural Networks(PINNs) \cite{Raissi19}, which have been applied in multiple contexts, including  solid mechanics \cite{Raissi21}. 
Other DL-based ROM techniques have come as an inspiration to handle the complex reduction process of dynamical systems, unveiling low-dimensional features from black-box data streams \cite{Guo2018,Guo2019,Carlberg20}. 
A peculiar perspective has been taken by 
Fresca et al.~\cite{FrescaManzoni21,fresca2020POD}, 
who proposed the non-intrusive POD-DL-ROM technique to address the fast simulation of parametrized differential problems. POD-DL-ROMs suitably combine a preliminary POD projection of snapshots, a Convolutional Autoencoder (AE) and a Deep Feedforward Neural Network (DFNN) to enable the construction of an efficient ROM, whose dimension is as close as possible to the number of parameters upon which the solution of the differential problem depends. The encoder part performs an operation of feature extraction forcing the high dimensional data to be reduced, at the {\em bottleneck} layer, to few reduced variables. We highlight that this approach builds on the idea that the system dynamics develops on a low dimensional invariant manifold, setting a clear parallel with the DPIM. 
Compared to other surrogate models exploiting machine/deep learning algorithms, a distinguishing feature of DL-ROMs is their capability to compute, at testing time, the whole solution field, for any new parameter instance and time instant, thus enabling the extremely efficient evaluation of any output quantity of interest depending on the solution field, and generating a truly virtual twin of the structure analysed.

In \cite{DLijnme}, an extension of the DL-ROM technique has been proposed and applied to microstructures, showing very good predictive capabilities. In particular, a first dimensionality reduction of the data is achieved 
by means of a POD-Galerkin (POD-G) ROM in order to reduce the costly FOM data generation phase, thus defining the POD-G DL-ROM technique.
An alternative approach is represented by the SINDy method (Sparse Identification of Dynamical Systems), a regression technique for extracting dynamics from time-series data  \cite{kaiser2018sparse,brunton2016discovering, brunton2016sparse}. In the past few years, SINDy has been widely applied to identify models of fluid flows, convection phenomena, structural modelling, and many others.
The main appeal of the SINDy approach is to generate, by combining a set of pre-defined (analytical) functions collected in a suitable dictionary, an 
explicit ROM in the form of first-order differential equations, which can be then manipulated with standard numerical tools. Even if applications have been limited so far to rather small problems, very recent works 
have fostered the use of SINDy in combination with AE neural networks \cite{champion2019data} and preliminary POD reduction.
The SINDy approach, with these recent extensions, and the aforementioned POD-DL-ROM, are thus converging to a similar framework, the difference being the 
way in which the reduced dynamics is integrated and approximated.
Convolutional or deep AEs thus emerge as distinctive and powerful tools for dimensionality reduction of dynamical systems.

However, a certain sense of distrust still pervades the community of nonlinear dynamics, and computational mechanics in general, towards these emerging DL approaches which are seen as black boxes lacking reliable foundations. Nevertheless, the strict connection with the DPIM emerges as quite evident, as discussed for instance in \cite{Chatzi21,yang21}. In the latter contribution, a data-driven identification of NNMs via
physics-integrated DL is proposed; in this case, the architecture integrates prior physics knowledge of the
NNMs by embedding physics-based constraints.
Intrinsic coordinates are required to have modal properties of a dynamical system, including statistical uncorrelations between desired modal displacements and modal velocities, and potential sparsity of dynamics in the intrinsic space. Their ideas are tested on systems
having 2 and 3 DOFs.

The POD-DL-ROM can be considered as the data-driven counterpart of the DPIM. Even if the two techniques arise from different perspectives, they share several analogies, thus generating a sort of duality. Indeed, they both rely on the use of a Full Order Model (FOM) -- typically represented by the FEM -- 
and they try to extract salient features of the problem being analyzed by reducing the model to an intrinsic space with very limited dimension. The reduced variables govern the evolution of the dynamics on invariant manifolds.
This reduction (or encoding) phase is performed offline with dedicated software and hardware as it can be memory and CPU intensive. 
A second shared phase is the integration of the reduced dynamics, which is fast and can be carried out almost in real-time. While DPIM achieves this latter task by integrating the reduced differential equations, e.g.\ with continuation techniques, POD-DL-ROMs adopts a DFNN, which directly approximates the parameter-reduced solution map. Finally, the whole field of FEM nodal displacements of the original FOM can be reconstructed (decoded) with great accuracy.  
For instance, the ROMs generated with both approaches investigated herein can be queried for 
the stress field at locations that do not need to be specified a priori, provided the query deals with parameter values falling within the range of (or not too far from) the training parameter set. 

Some features however differentiate the two approaches. The DPIM is an almost exact procedure which, for a given FOM with fixed parameters and
forcing type, generates a ROM which is highly predictive for all initial conditions and
forcing levels. However, it does not generalise automatically to other physics and
nonlinearities. 
Conversely, the DL-based ROM accommodates a whole range of parameters selected by the user in the training phase, can be extended almost automatically to coupled physics, and has real-time performances. However, it still suffers from limited predictive ability and does not guarantee the same level of accuracy as the DPIM. Indeed, such a ROM strategy account for the problem physics only through the snapshots data, employed to train the neural network.

As a natural consequence of these remarks, the aim of this contribution is twofold.
After a short introduction to the governing equations in Section \ref{sec:formulation} and of the DL-based ROM method in Section \ref{sec:ROM}, we further improve the POD-DL-ROM technique \cite{FrescaManzoni20,FrescaManzoni21} in Section \ref{sec:abscissa} through the use of an arc-length abscissa defined on the Frequency Response Functions (FRFs) to enable the analysis of complex interaction problems in real nonlinear microstructures 
showing, e.g., internal resonances and coupled electromechanical interactions inducing autoparametric effects.
In parallel, we explore in-depth
the relationship between the DL-ROM and the DPIM approaches in Section \ref{sec:dpim}. We devote particular attention
to the convergence of the POD-DL-ROM with respect to the dimension of the reduced space, analyzing it in terms of both the FRF and the Invariant Manifolds predicted. 
The applications are proposed in Section \ref{sec:applications}.
In the first two examples addressed, which focus on mechanical geometric nonlinearities, the DPIM simulations are used as a reference solution due to their known accuracy, and the convergence of the POD-DL-ROM approximation to the DPIM one is investigated. In the final example, on the contrary, in order to show that the POD-DL-ROM easily generalises to complex multi-physics without losing its non-intrusive nature, we exploit a commercial code to generate snapshots of an electromechanical gyroscope, aiming at simulating its oscillating behavior in real-time.

\section{Problem formulation}
\label{sec:formulation}

We will focus on mechanical structures subjected to periodic forcing 
that undergo large transformations but still experience only small strains,  a condition which is well described by the Saint Venant-Kirchhoff constitutive model. As detailed in \cite{DLijnme}, the space discretization of the governing equations by means of finite elements yields a system of coupled nonlinear differential equations representing the Full Order Model (FOM):
\begin{subequations} 
	\label{eq:PPV_d2}
	\begin{align}
		&\bfM \ddot{\bfU}(t) + \bfC\dot{\bfU}(t)  + \bfK\bfU(t) + \bfG(\bfU,\bfU) + \bfH(\bfU,\bfU,\bfU) = 
		\bfF(t;\bfmu), \qquad t \in (0,T),\\
		&\bfU(0)=\bfU(T), \qquad
		\dot{\bfU}(0)=\dot{\bfU}(T), 
	\end{align}
\end{subequations}
where the vector $\bfU(t) \in \mathbb{R}^{N_h}$ collects the $N_h$ unknown displacements nodal values; 
$\bfM \in \mathbb{R}^{N_h \times N_h}$ is the mass matrix; $\bfC=(\omega_0/Q)\bfM$ is the Rayleigh mass-proportional damping matrix with $\omega_0$ reference eigenfrequency and $Q$ quality factor;
$\bfF(t;\bfmu) \in \mathbb{R}^{N_h}$ is the nodal force vector which depends on the vector of $n_{\bfmu}$ parameters $\bfmu \in \cP \subset \Rbb^{n_{\bfmu}}$ and expresses the actuation due, in general, to a multiphysics coupling e.g.\ with piezolectricity or electrostatics. In particular, the vector $\bfF(t;\bfmu)$ depends on the angular frequency $\omega$ 
of the actuation and is nonlinearly modulated in amplitude by the coefficient $\beta$.

The internal force vector has been exactly decomposed in linear, quadratic, and cubic power terms of the displacement: $\bfK \in \mathbb{R}^{N_h\times N_h} $ is the stiffness  matrix related to the linearized system, while $\bfG \in \mathbb{R}^{N_h}$ and $\bfH \in \mathbb{R}^{N_h}$ are  vectors given by monomials of second and third order, respectively. 
We stress that the components of these vectors can be expressed using an indicial notation as
\[
G_i=\sum_{j,k=1}^{N_h}G_{ijk} U_j U_k, \quad   H_i=\sum_{j,k,l=1}^{N_h}H_{ijkl}U_j U_k U_l, \qquad i=1,\ldots, N_h.
\] 
Eqs.~\eqref{eq:PPV_d2} represent our high-fidelity FOM which depends on the input parameters $\bfmu$. Our goal is the efficient numerical approximation of the solution manifold:
\begin{equation}
	\label{eq:sethifi}
	\cS=\left\{\bfU(t;\bfmu) \, : \, t \in [0,T)\, ,  \, \bfmu \in \cP \subset \Rbb^{n_{\bfmu}} \right\}\subset \Rbb^{N_h}.
\end{equation}
The numerical solution of Eqs.~\eqref{eq:PPV_d2} to compute the steady state response is a challenge by itself for large scale problems. One option is the use of time marching methods to simulate a sufficiently large number $N_c$ of cycles, where $N_c$ is typically inversely proportional to the damping. 
This technique resorts to robust algorithms implemented in most of the commercial software, but when damping is very small, as in most MEMS devices, the computational effort may be unaffordable. 
In other approaches, like the Harmonic Balance (HB) method \cite{actuators21,KerschenHB2015}, the 
unknown displacements are expressed as the sum of Fourier components thus automatically respecting periodicity conditions. However, their implementation requires dedicated codes and non-standard computing facilities.

\section{Data-driven Reduced Order Modeling through DL-ROMs}
\label{sec:ROM}

In this section, we briefly review the construction of deep learning-based ROM (DL-ROM) technique and its extension to the POD-enhanced version (POD-DL-ROM).

The DL-ROM technique \cite{FrescaManzoni21} is extremely efficient and is able to model highly nonlinear problems by identifying the manifold underlying the dynamics of the system in a complete data-driven and black-box, non-intrusive way. On the other hand, its data-driven nature implies that the sampled data provided must span the parameter space of interest and contain all the information necessary to accurately approximate the solution manifold. In this way, the size of the training dataset increases with the number of parameters considered in the FOM. 

DL-ROMs aim at approximating the map $(t, \bfmu) \rightarrow {\bf U}(t, \bfmu)$ by describing both the trial manifold, approximating (\ref{eq:sethifi}), and the reduced dynamics through deep neural networks, which are trained on a set of FOM snapshots. In particular:
\begin{itemize}
\item to map the FOM solutions in a low-dimensional coordinates vector representation (\textit{encoding stage}), we use the encoder function of a Convolutive AutoEncoder (CAE)
\begin{equation}
\label{eq:encoder}
 {\mathbf{\tilde{z}}_n}(t; \boldsymbol{\mu}, \boldsymbol{\theta}_{E}) = {\mathbf{f}}_{n}^E(\mathbf{U}(t; \boldsymbol{\mu}); \boldsymbol{\theta}_{E});
\end{equation}
where $\mathbf{f}_{n}^E(\cdot; \boldsymbol{\theta}_E) : \mathbb{R}^{N_h} \rightarrow \mathbb{R}^n$ is obtained as   the composition of several layers (some of which are convolutional), depending upon a vector ${\boldsymbol{\theta}}_E$ of parameters;
\item to describe the system dynamics (\textit{reduced dynamics learning}), the intrinsic coordinates of the ROM approximation are defined as
\begin{equation}
\label{eq:reduced_dynamics}
{\mathbf{z}}_n(t; \boldsymbol \mu, \boldsymbol{\theta}_{DF}) = {\boldsymbol{\phi}}_n^{DF}(t; \boldsymbol \mu, \boldsymbol{\theta}_{DF}), 
\end{equation}
where ${\boldsymbol{\phi}}_n^{DF}(\cdot; \cdot, \boldsymbol{\theta}_{DF}) : {\mathbb{R}}^{(n_{\mu} +1)} \rightarrow {\mathbb{R}}^n$ is a DFNN, consisting in the subsequent composition of a nonlinear activation function,  applied to a linear transformation of the input, multiple times \cite{goodfellow2016deep}.  Here $\boldsymbol{\theta}_{DF}$ denotes the vector of {parameters} of the DFNN, collecting all the weights and biases of each layer of the DFNN; 
\item  to model the reduced nonlinear trial manifold  ${\mathcal{S}}_h^n \approx \mathcal{S}$ (\textit{decoding stage}) we employ the decoder function of a CAE  \cite{lecun1998gradient, hinton1994autoencoders}, that is, 
\begin{equation}
\label{eq:nonlinear_manifold}
\begin{split}
\tilde{\mathcal{S}}_h^n = \{ {\mathbf{f}}^D_h(\mathbf{z}_n(t; \boldsymbol{\mu}, {\boldsymbol{\theta}_{DF}}); \boldsymbol{\theta}_{D})
\, : \,  \mathbf{z}_n(t; \boldsymbol{\mu}, \boldsymbol{\theta}_{DF}) \in {\mathbb{R}}^{n} \, , \,   t \in [0, T) \, , \, \boldsymbol{\mu} \in \mathcal{P} \} \subset \mathbb{R}^{N_h},
\end{split}
\end{equation}
where ${\mathbf{f}}^D(\cdot; {\boldsymbol{\theta}}_D) : {\mathbb{R}}^n \rightarrow {\mathbb{R}^{N_h}}$ depends upon a vector ${\boldsymbol{\theta}}_{D}$ collecting all the corresponding weights and biases. 
\end{itemize}

The DL-ROM approximation $\mathbf{\tilde{U}}(t; \boldsymbol \mu) \approx \mathbf{U}(t; \boldsymbol \mu)$  is then given by
\begin{equation}
\mathbf{\tilde{U}}(t; \boldsymbol \mu, \theta_{DF}, \theta_D) = {\mathbf{f}}_h^D({\boldsymbol{\phi}}_n^{DF}(t; \boldsymbol{\mu}, {{\boldsymbol{\theta}}_{DF}}); \boldsymbol{\theta}_{D}),
\label{eq:reconstructed_solution}
\end{equation}
and is computing by solving a suitable optimization problem, in the  variable $\boldsymbol{\theta} = (\boldsymbol{\theta}_{E}, \boldsymbol{\theta}_{DF}, \boldsymbol{\theta}_{D})$ (see, e.g., \cite{FrescaManzoni21} for further details).

The POD-DL-ROM technique \cite{fresca2020POD} consists in applying the DL-ROM technique to the intrinsic coordinates of a linear trial manifold generated through randomized singular value decomposition (rSVD) and approximating $\mathcal{S}$; In particular, to reduce the dimensionality of the snapshots and avoid feeding training data of very large dimension $N_h$, POD is first applied -- realized through randomized SVD (rSVD) \cite{halko2011finding} -- to the snapshot set; then, a DL-ROM is built to approximate the map between $(t, \bfmu)$ and the POD generalized coordinates $\bfu_N(t; \bfmu) = \mathbf{V}_N^T \mathbf{U}(t; \bfmu)$. 

\begin{figure}[htb]
	\centering
	\includegraphics[width =.85\textwidth]{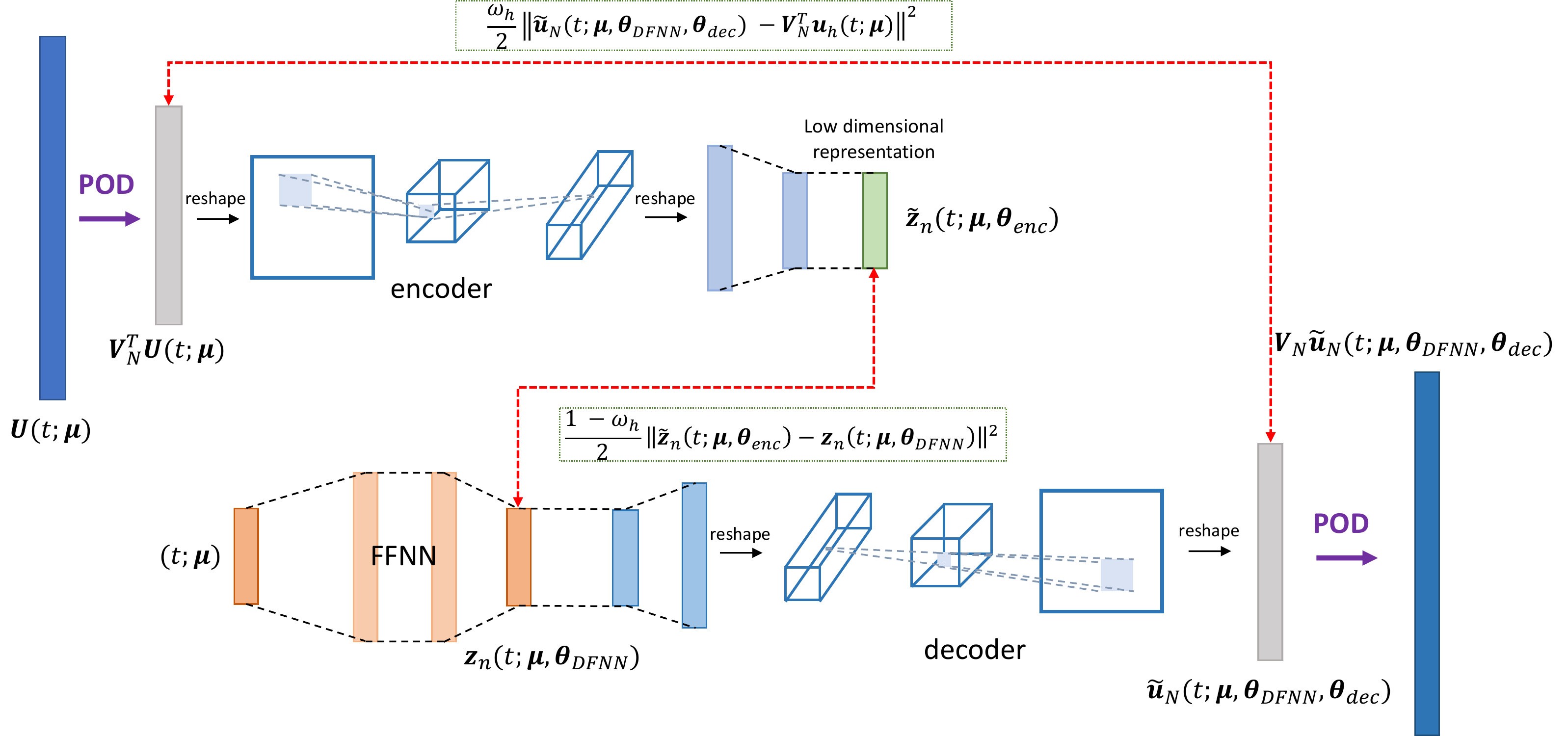}
	\caption{POD-DL-ROM architecture. Starting from the FOM solution $\mathbf{U}(t; \boldsymbol{\mu})$, the intrinsic coordinates $\mathbf{V}_N^T \mathbf{u}_h(t; \boldsymbol{\mu})$ are computed by means of rSVD;  the neural network provides their approximation $\mathbf{\tilde{u}}_N(t; \boldsymbol{\mu})$ as output. The reconstructed solution $\mathbf{\tilde{u}}_h(t; \boldsymbol{\mu})$ is then recovered through the basis matrix.}
	\label{fig:POD-G DLROM}
\end{figure}

Note that, in the POD-DL-ROM framework, the encoding stage results from the projection of FOM snapshots onto the linear trial manifold defined by $\mathbf{V}_N$ and the compression performed by the encoder function of the CAE. Similarly, the decoding stage is formed by the application of the decoder function and the multiplication for the POD basis matrix. 
The architecture of the POD-DL-ROM neural network, employed at training time, is the one shown in \figurename~\ref{fig:POD-G DLROM}; note that,  at testing time, as in the DL-ROM technique, we can discard the encoder function.

A variant of the POD-DL-ROM approach, called POD-G DL-ROM, has been recently benchmarked on several MEMS structures including beams, arches and mirrors in \cite{POD21} with the aim of reducing the cost of the training phase. A limited number of FOM analyses are used to identify an optimal linear basis via SVD, and the ROM obtained projecting the governing equation on this basis is then employed to train the neural networks.
However, to assess the capacity of the POD-DL-ROM 
to tackle other, more involved physical problems, the initial formulation is adopted in this work.

\subsection{Arc-length on the FRF as control parameter}
\label{sec:abscissa}

In applications to microstructures, one of the most meaningful outputs of the analysis is the Frequency Response Function (FRF) which defines the steady-state periodic response of the dynamical system. In a FRF, a selected output like, e.g., the maximum midspan deflection of a beam, or the rotation amplitude of a micromirror, is plotted as a function of the actuation intensity and frequency. In linear systems, the FRF is a single-valued function of the forcing frequency $\omega$ that is hence utilized as the ordering parameter also for the snapshot matrix. When nonlinearities are present, this property is generally lost. For instance, in systems behaving as simple duffing oscillators with hardening (or softening) properties, the phase of the response with respect to the forcing signal has been exploited in \cite{DLijnme} as order parameter.

\begin{figure}[t!]
	\centering
	\includegraphics[width = .85\linewidth]{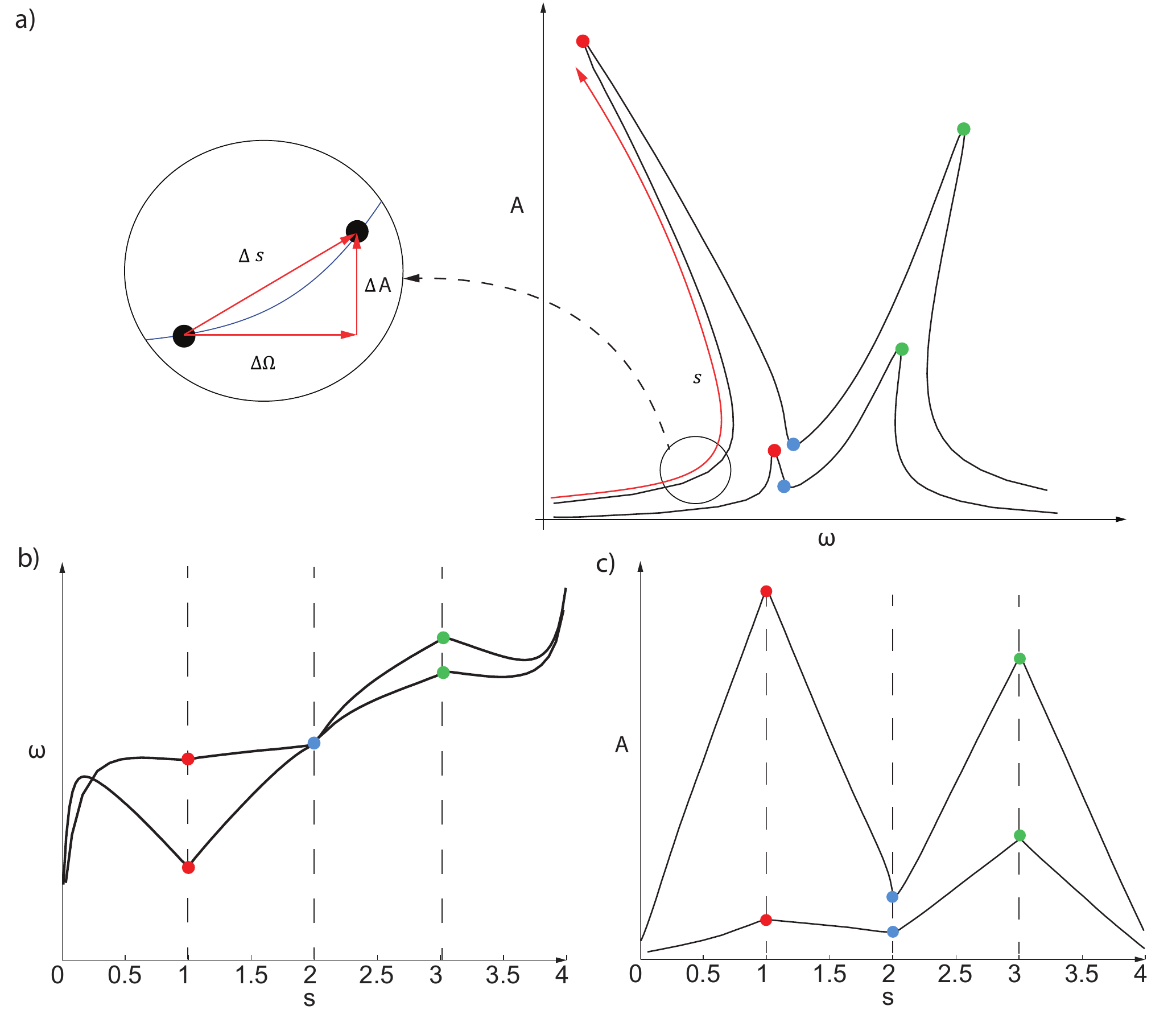}
	\caption{Arch-length abscissa. Figure a): the arch-length abscissa is first computed along each FRF. Figure b) and c): the new $\omega(s)$ and $A(s)$ are single valued functions. Peaks and valleys of the FRF are aligned by rescaling the arc-length abscissa in four different regions.}
	\label{fig:arclength}
\end{figure}

In more involved applications, even the phase does not ensure uniqueness of the response, and a more general approach is needed. Inspired by the continuation technique used to integrate the ROM in the DPIM approaches, we opt for the use of the arc-length abscissa along the FRFs.
Let us consider for instance the typical FRF for a 1:2 internal resonance of a shallow arch, as illustrated in Figure \ref{fig:arclength}a), where the curves
for two different forcing levels are plotted.
In order to express both the frequency and the amplitude as single valued functions, we first compute the arc-length abscissa along the FRF. However, the total length of the FRFs strongly depends on the forcing levels; moreover, the introduction of this abscissa induces a {\it misalignment} of physically relevant phenomena, like e.g. resonances and mode couplings.
These are, on the contrary, naturally clustered at very similar values of the frequency and phase with positive fall-outs in the DL-ROM.

As a solution to this issue, we introduce a piecewise normalization of the arc-length abscissa that synchronizes maxima and minima of the FRFs. 
Four regions are identified, as illustrated in Figure \ref{fig:arclength}, and the arc-length of each region is rescaled from 0 to 1, yielding a total arc-length of 4.  
Even if this procedure requires some insight into the physics of the response, it is however very general and powerful, thus representing a major improvement with respect to other procedures which would not be applicable in this context \cite{DLijnme}.

We stress moreover that the arc-length abscissa is not only used as control parameter
in the snapshots passed to the DL-ROM; indeed, it also
easily allows to retrieve, a posteriori, the FRF frequency-amplitude relationship through interpolation on the training data. 

\FloatBarrier
\section{Direct Parametrization of Invariant Manifolds}
\label{sec:dpim}

The parametrization method of invariant manifolds, introduced by Cabr{\'e}, Fontich and de la Llave in~\cite{Cabre1,Cabre2,Cabre3}, has been very recently applied to large FEM models of nonlinear vibrating structures \cite{vizza2020,NNM21,opreni22}. In particular, in \cite{opreni22} it is shown that lightly damped and mildly forced systems can be parametrized 
with very good accuracy by
considering only the autonomous part of the equations and adding a posteriori the forcing
by projecting $\bfF$ in Eq.\eqref{eq:PPV_d2} on the linear master mode. 
This is the approach retained in the present work.
The general idea is to reduce the dynamics to an invariant manifold tangent to the eigenvectors selected as master modes. Let us assume that $n$ master coordinates are selected, with $n\ll N_h$. These master coordinates are linked to their corresponding vibration modes and the searched invariant manifold is the nonlinear continuation of the subspace spanned by the $n$ vibration modes. 
Working in the phase space in terms of both displacements and velocities, the invariant manifold becomes $2n$-dimensional; to describe the reduced dynamics on this manifold, $2n$ \textit{normal} coordinates $\bfz$ are introduced. The $2N$ original coordinates $\bfU$ (nodal displacements) and $\bfV$ (nodal velocities) are 
then expressed as a function of the new normal coordinates $\bfz$ as
\begin{align}
\label{eq:mappings_compact}
\bfU = \bfPsi(\bfz),
\quad
\bfV = \bfUps(\bfz),
\end{align}
where the two nonlinear mapping functions $\bfPsi(\bfz)$ and $\bfUps(\bfz)$ are time-independent unknowns to be computed. 
The reduced dynamics of the autonomous system governs the evolution onto the corresponding invariant manifold. This latter is also unknown at this stage, and can be expressed as:
\begin{equation}
\label{eq:reddyn_compact}
\dot{\bfz} = \bff(\bfz).
\end{equation}
The aim of the method is to compute $\bfPsi$, $\bfUps$ and $\bff$. 
The procedure consists in differentiating Eq.~\eqref{eq:mappings_compact} with respect to time
\begin{align}
\label{eq:diff}
\dot{\bfU} 
= \nabla_\bfz\bfPsi(\bfz) \, \dot{\bfz} 
= \nabla_\bfz\bfPsi(\bfz) \,  \bff(\bfz),
\qquad
\dot{\bfV} 
= \nabla_\bfz\bfUps(\bfz) \, \dot{\bfz} 
= \nabla_\bfz\bfUps(\bfz) \,  \bff(\bfz),
\end{align}
and substituting in the equations of motion \eqref{eq:PPV_d2} -- suitably rewritten as first-order system -- to obtain the so called {\it invariance equations} 
\begin{subequations}
\begin{align}
& \bfM \;\nabla_\bfz \bfUps(\bfz)\; \bff(\bfz) + \bfC \bfUps(\bfz) + \bfK \bfPsi(\bfz)
 +\bfG(\bfPsi(\bfz),\bfPsi(\bfz)) + \bfH(\bfPsi(\bfz),\bfPsi(\bfz),\bfPsi(\bfz)) = \bfm{0},
\\
& \bfM\nabla_\bfz\bfPsi(\bfz)\; \bff(\bfz) = \bfM\bfUps(\bfz).
\end{align}\label{eq:invariance_compact}
\end{subequations}
These nonlinear equations can be solved using asymptotic expansions in the unknowns (the reduced coordinates $\bfz$), as proposed in~\cite{Cabre3,Haro}. The procedure for large FEM models
has been detailed in \cite{vizza2020,NNM21}, extended to high-order expansions in \cite{vizzaccaro2021high} and
to forced systems in \cite{opreni22}, where it has been extensively validated against FOM HB solutions. An open-source version of the code has been published in \cite{Morfe} and is now a reference for geometrical nonlinearities with an impressive list of successful applications ranging from internal resonance, 
autoparametric excitation and folding of the manifolds due to extremely large transformations, just to mention a few. However, the DPIM formulation has not been extended yet to coupled multiphysics.

One of the greatest benefits of the DPIM is the generation of a ROM in the form of an explicit, small-size system of first-order differential equations -- see Eq.\eqref{eq:reddyn_compact} --
which can be solved with different approaches.
MEMS structures, and in general periodically excited systems, are characterized by their steady state regime. Since damping is typically very low, time-marching methods are less appealing than Harmonic Balance, Collocation or Shooting techniques that allow computing the steady state response directly. Furthermore, these methods are compatible with continuation approaches that allow reconstructing the whole FRF, with both stable and unstable branches. 
Several packages suitable for small-scale problems like the one given by DPIM are freely available. Among them we cite \texttt{Auto07p} \cite{Doedel}, \texttt{Manlab} \cite{Guillot1,Guillot2}; \texttt{Nlvib} \cite{NVLIB}, \texttt{MATCONT} \cite{dhooge2006matcont},
\texttt{COCO}  \cite{COCO} and  \texttt{BifurcationKit} \cite{BifurcationKit}, an emerging toolkit for Julia language. All the continuation solutions discussed in this work have been obtained using \texttt{MATCONT}.

\subsection{Duality with the DL-ROM in a simple case}
\label{sec:DPIM}

\begin{figure}[h]
	\centering
	\includegraphics[width = .8\linewidth]{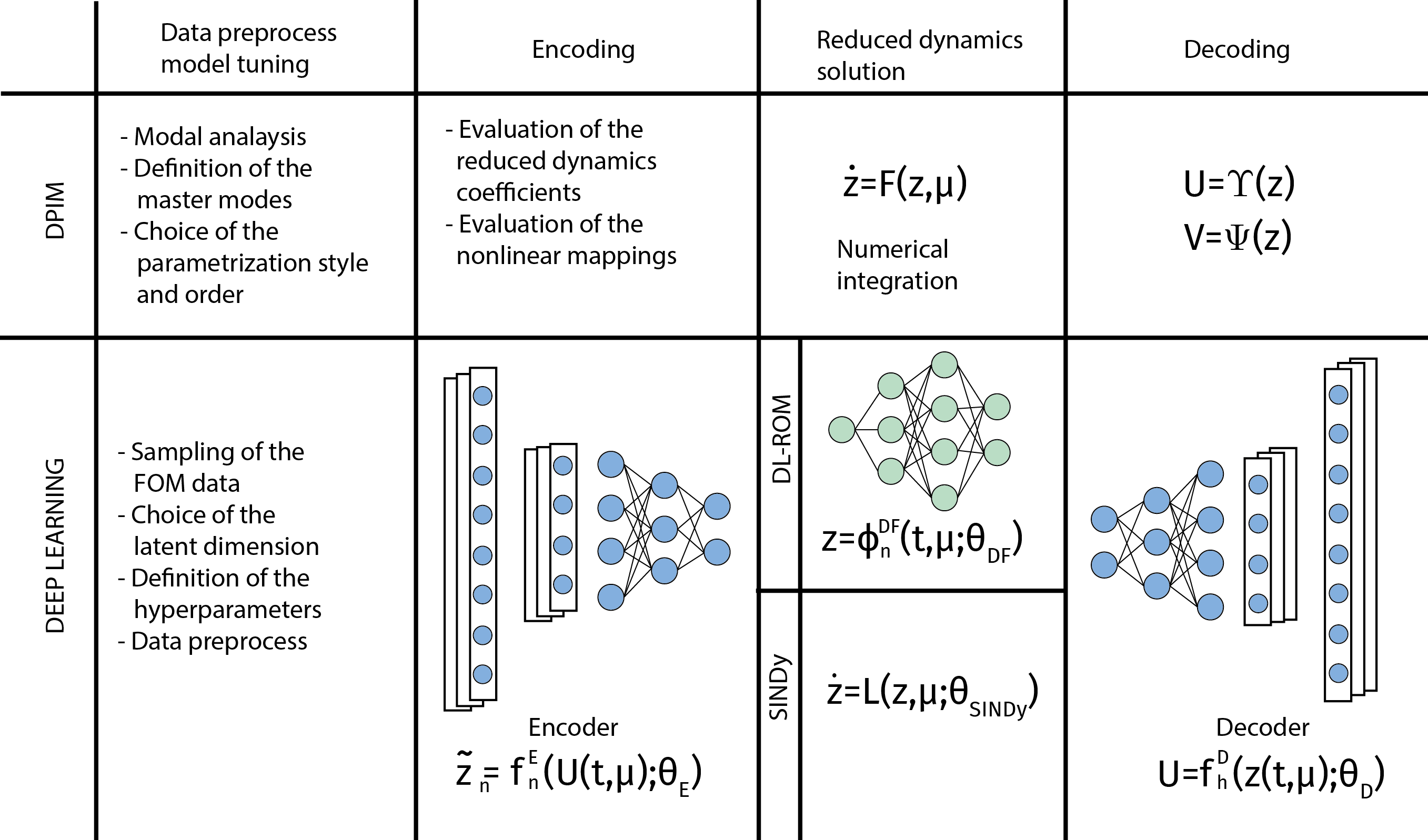}
	\caption{
	Schematic comparison between DPIM and DL methods. 
	One can identify a pre-processing stage, an encoding stage, a reduced dynamics solution stage and, finally, a decoding stage. For DL methods, in this scheme we distinguish between black box approaches for the reduced order dynamics, as in the DL-ROM \cite{FrescaManzoni21}, and model based approaches, as in SINDy proposed by Brunton et al. \cite{brunton2016discovering} 
	}
	\label{fig:DPIMDL}
\end{figure}

In order to better highlight the striking duality/analogy between the DL-ROM and the DPIM, we now focus
on the real formulation of the DPIM proposed in \cite{vizza2020}, and restrict the attention to a single mode reduction in a mechanical system
with asymptotic expansion truncated at second order. 
The basic ingredients of the procedure can be expressed in terms of the real coordinates $R$ and $S$, i.e.\ the reduced {\it normal coordinates}, which are the counterpart of $\bfz$ in Eqs.\eqref{eq:mappings_compact}
and \eqref{eq:reddyn_compact}. In this case, the nonlinear change of variables reads:
\begin{align}
\label{eq:reconstr1}
\bfU & = \bfPhi R + \bfa R^{2} + \bfb S^{2}+\bfc RS,
\\
\label{eq:reconstr2}
\bfV & = \bfPhi S + \bfalpha R^{2} + \bfbeta S^{2} +\bfgamma RS, 
\end{align}
and allows expressing the FOM nodal values for displacements $\bfU$ and velocities $\bfV$ from the normal coordinates. The reduced model can be formulated in terms of two first-order equations as follows:
\begin{align}
\label{eq:rom1}
\dot{R} & =S,
\\
\label{eq:rom2}
\dot{S} & = -\frac{\omega_0}{Q} S -\omega_0^2R - A R^3 - B RS^2
- C R^2 S + \beta \cos(\omega t),
\end{align}
All the vectors and coefficients in Eqs.\eqref{eq:reconstr1}-\eqref{eq:rom2}
are computed
in the parametrization procedure while $\beta$ denotes the forcing level. 
Let us now underline some contact points between the DL-ROM and the DPIM (see Figure \ref{fig:DPIMDL}): 
\begin{enumerate}
    \item   All the coefficients and vectors in Eqs.\eqref{eq:rom1}-\eqref{eq:rom2}
are computed offline starting from the FOM, in a preliminary phase, that can be seen as the equivalent of the {\it encoding} phase in the DL-ROM.
The costly offline training can be performed on dedicated platforms and software. Both approaches hence are based on a FOM, which is typically built exploiting a finite element discretization.
\item The ROM Eqs.\eqref{eq:rom1}-\eqref{eq:rom2} emerges as the counterpart of the DFNN in Eq.\eqref{eq:reduced_dynamics} and can be integrated online with almost real-time performance when the model is queried with specific values of the forcing parameters.
\item Eqs.\eqref{eq:reconstr1}-\eqref{eq:reconstr2} represent the parallel of the {\it decoding phase} which reconstructs the global nodal fields 
starting from 
the online integration of the ROM.
Thus both ROMs can reproduce the same richness in details of the original FOM,
since the decoding phase allows to generate a full field information.
\end{enumerate}
It is worth recalling that the SINDy approach \cite{kaiser2018sparse,brunton2016discovering, brunton2016sparse} has been recently coupled with autoencoders for order reduction \cite{champion2019data} and shares several of the similarities in common between the DPIM and the DL-ROM. However, differently from the latter, SINDy expresses the reduced dynamics in the form of an ODE system whose coefficients are estimated by means of sparse identification procedures.

\FloatBarrier
\section{Applications}
\label{sec:applications}

In this section, we discuss three different applications: a micromirror, a shallow arch showing 1:2 internal resonance, and an electromechanical gyroscope.
In the first two examples addressed, 
which focus on geometrical nonlinearities, the DPIM simulations are used as reference solutions, and the convergence of the DL-ROM to the DPIM is investigated. A detailed discussion on the performance of the DPIM 
in these kind of applications has been developed
in \cite{NNM21,vizzaccaro2021high} through extensive validations against the HB approach.

In the final example, on the contrary, we exploit a commercial code to generate FOM snapshots 
of an electromechanical gyroscope in order to show that the construction of a DL-ROM can be easily generalized to the case of complex multi-physic problems, still retaining its non-intrusive peculiar character. In this third example, validation will be performed against 
the results of the commercial FOM itself.

\paragraph{Reconstruction of the whole response}

As highlighted in the previous Sections, the methods analysed are not limited to {the evaluation of} the FRFs of selected quantities like, e.g., the rotation angle or the midspan displacement, but they accurately 
regenerate, through the decoding phase, an approximation $\tilde{\bfU}$ of the whole displacement field of the original FOM. 
Next, $\tilde{\bfU}$ can be post-processed to generate all the desired outputs, a typical example being the monitoring the evolution of stresses  in the structure to check the admissibility of the given design \cite{DLijnme}.
In our case, we will use $\tilde{\bfU}$ to reconstruct the FRFs for the master and the slave modes,
in order to investigate the convergence of the DL-ROM to the DPIM results.

Given the estimate of the nodal displacement vector $\tilde{\bfU}$  at any instant
$t$, we define a generalized modal coordinate $u_{i}(t)$ 
and modal amplitude $A_i$ according to:
\begin{align}
\label{eq:umod}
 u_{i}(t)=\bfPhi^T_i \bfM  \tilde{\bfU}(t)
 \qquad
A_i=\max_{t\in [0-T]}(|u_i(t)|)
\end{align}
in which $\bfM$ is the mass matrix, $\bfPhi_i$ is the i-th eigenmode and $T$ is the period. Modal coordinates are also used to generate error measures for each mode:
\begin{align}
\label{eq:err1}
E^r_{i}=\|\uP_i(t)-\uDL_i(t)||_{L_2}/||\uP_i(t)||_{L_2}
\\
\label{eq:err2}
E_{i}=\|\uP_i(t)-\uDL_i(t)||_{L_2}/||\uP(t)||_{L_2}
\end{align}
where $\uP_i$ denotes the DPIM solution, $\uDL_i$ the DL-ROM one,
and $(\uP)^2=\tilde{\bfU}^T\bfM\tilde{\bfU}$.
The proposed errors are essentially time averages of the instantaneous error of the DL-ROM with respect to the DPIM reference. The former one, $E^r_{i}$, is a relative indicator for a given mode, but it does not account for the absolute importance of the mode itself
in the global response. The latter one, $E_{i}$, instead normalizes the error with respect to a global displacement measure.

\paragraph{Computation of manifolds}

The underlying key idea and assumption of the two reduction techniques addressed in this contribution is the existence of invariant manifolds on which the system 
response evolves as a function of the (few) master coordinates.  As a consequence, in all the examples we will analyse the existence of DL-ROM manifolds 
and their convergence, when applicable, to the DPIM ones. 
While plotting the manifold for the DPIM is straightforward as both displacements and velocities are directly accessible, see Eq.\eqref{eq:mappings_compact}, 
a DL-ROM provides, on the contrary, only a reconstruction of displacements.
In order to generate the manifold, the modal velocity $\vDL_i(t)$ can be computed resorting to a Fourier decomposition of the periodic $\uDL_i(t)$, and then differentiating it with respect to time. 
Indeed, this approach allows to smooth high-order components  through harmonic truncation, and is preferred, e.g., over finite differences that might generate noisy derivatives.

\subsection{Micromirror}
\label{sec:mirror}

Micromirrors represent one of the most promising new families of MEMS, with applications ranging from LIDARs \cite{Lidar} 
to pico-projectors or Augmented Reality (AR) lenses \cite{ARlenses}, e.g.\ Microsoft Hololens \cite{hololens}.

\begin{figure}[ht]
	\centering
	\includegraphics[width = .8\linewidth]{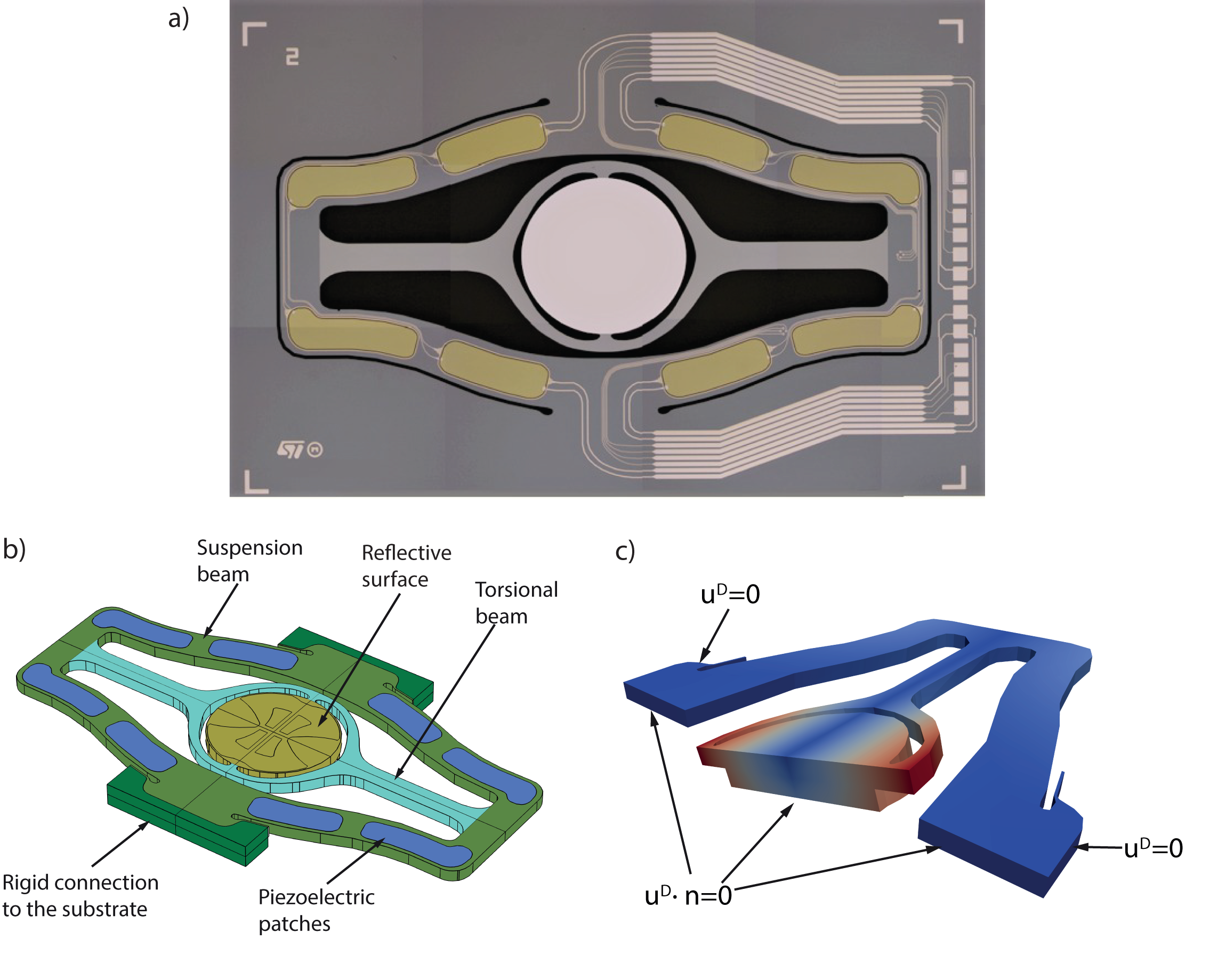}
	\caption{Micromirror: optical photo and schematic view}
	\label{fig:perseus}
\end{figure}
\begin{figure}[ht]
	\centering
	\includegraphics[width = .7\linewidth]{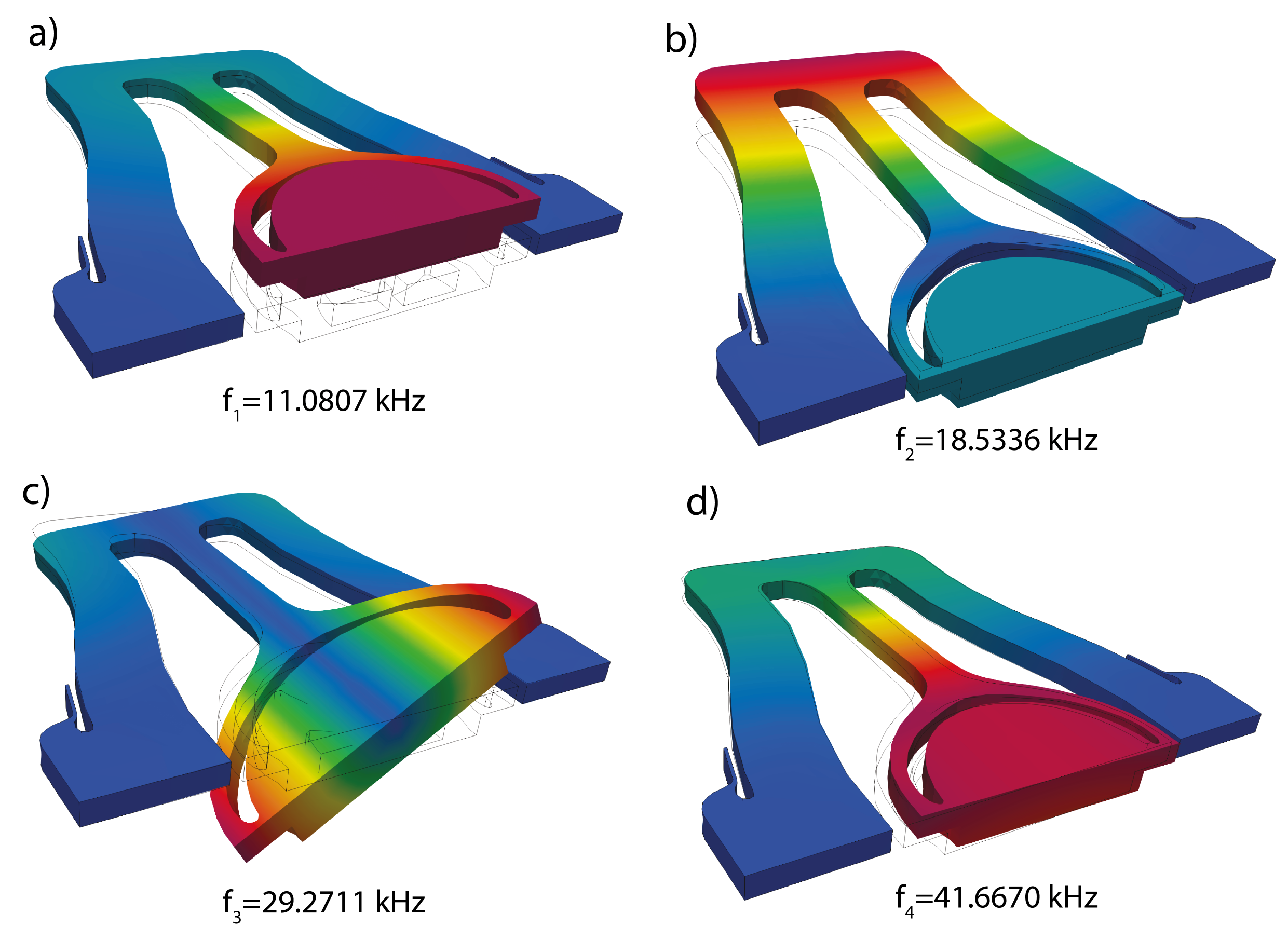}
	\caption{Micromirror. Figure a)-d) first four eigenmodes of the micromirror. Figure c): torsional master mode.}
	\label{fig:perseusmodes}
\end{figure}

Because of the inertial and geometrical effects triggered by large rotations, micromirrors are intrinsically nonlinear, and the prediction of their dynamic behaviour is essential to guarantee a proper design and control of the device during operations.
Recently, the authors have developed a large scale HB-approach in \cite{actuators21} for the analysis of piezo-actuated mirrors. However, the computational cost entailed by this approach limits its applicability during the design phase and for online monitoring, 
ultimately stimulating a frantic search for efficient ROMs. This example is a tough challenge not only for classical implicit condensation approaches \cite{ijnm19} that cannot deal with large rotations, but also
for the most advanced and recent techniques, like the DPIM approach. 
The torsional mode is not the lowest-frequency one and is not well separated from other eigenmodes, entailing the failure of the quadratic formulation of the DPIM utilized in \cite{NNM21}, therefore requiring a high-order expansion, as remarked by Vizzaccaro et al.~\cite{vizzaccaro2021high}. 
In this section, we will consider the MEMS micromirror illustrated in Figure~\ref{fig:perseus}, fabricated by ST Microelectronics.  The mirror is assumed to be made of isotropic polysilicon \cite{jmems04}, with density $\rho=2330$\,Kg/m$^3$, Young modulus $E=167$\,GPa and  Poisson coefficient $\nu=0.22$.
The mirror plate, a circle of diameter $1970$\micr, is suspended to a gimbal connected with a torsional beam along the rotation axis and two suspension beams on each side. Thanks to symmetry, only half of the device is modelled with the FEM, involving in total 9732 dofs;  Dirichlet boundary conditions are imposed as illustrated in Figure~\ref{fig:perseus}c) and on the symmetry plane. On the remaining boundaries, zero traction Neumann boundary conditions are instead imposed.
The mirror has been actuated with a fictitious body force proportional to the inertia forces of the master mode, 
$\bfF(t)=\bfM  \bfphi_3 \beta\cos(\omega t)$, in order to enable the comparison with the DPIM approach \cite{vizzaccaro2021high}. 
The first four eigenmodes are illustrated in Figure~\ref{fig:perseusmodes}, where also the eigenfrequencies obtained from a linear FEM eigenvalue analysis are listed. 
The torsional {\it master} mode, of interest herein, is the third one and has a frequency of $29 271$\,Hz.  A quality factor $Q=1000$ is considered in all the analyses.

\begin{figure}[h]
	\centering
	\includegraphics[width = .8\linewidth]{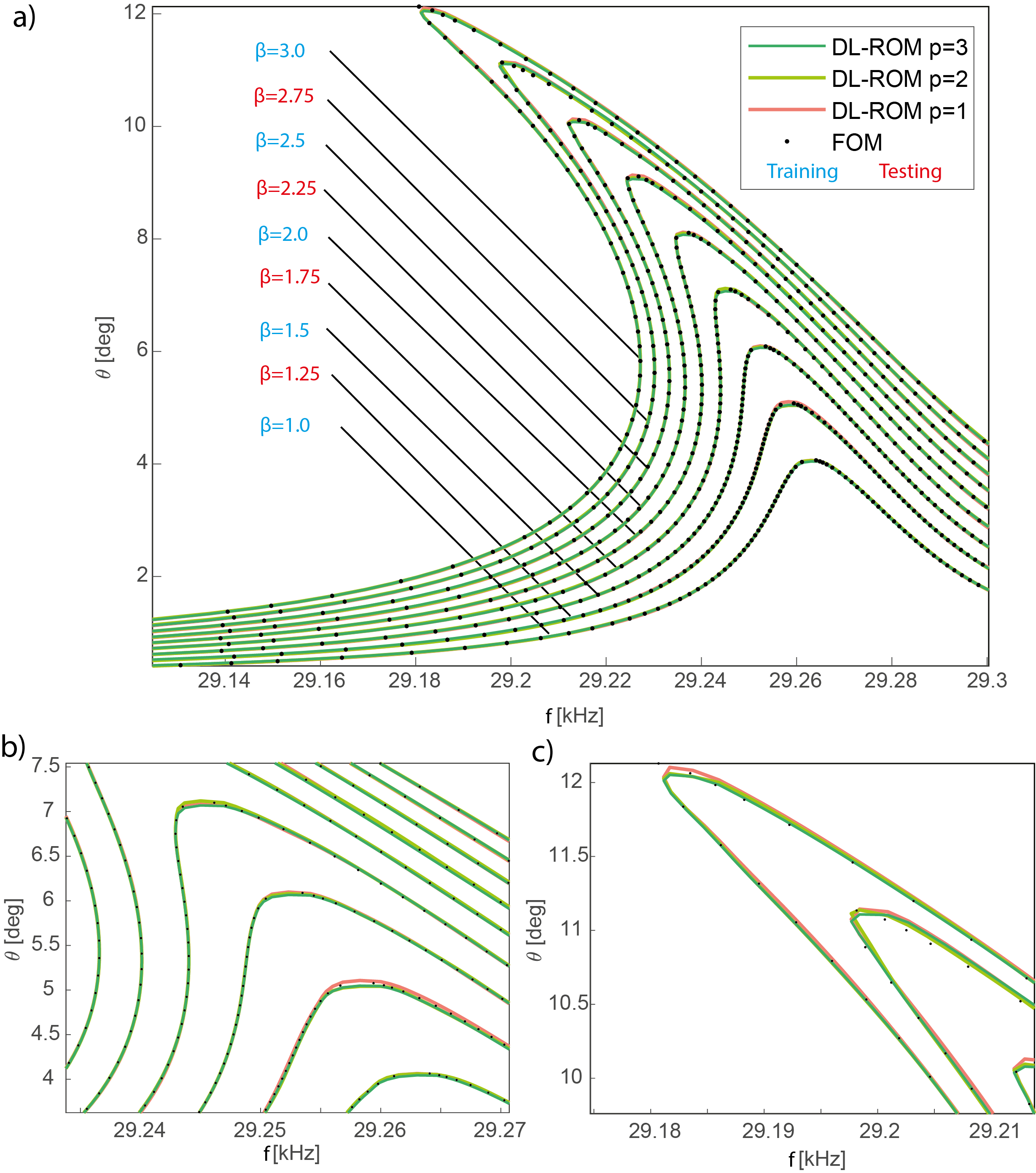}
	\caption{Micromirror. FRFs of the opening angle for different excitation levels $\beta$. Figure a): comparison between the FRFs obtained with the DPIM and the DL-ROM using 
	$p=1\shddot 3$. Figures b) and c): enlarged views of the peaks in the response. The training $\beta$ values are labelled in light blue while the testing ones are in red. 
	The FRFs here reported have been obtained with a testing dataset made of 112 700 instances. The inquiry of the DL-ROM is performed in less than 0.2 s using a Tesla V100 32GB GPU and an implementation in the Tensorlow DL framework.}
	\label{fig:perseusFRF}
\end{figure}

\begin{figure}[h]
	\centering
	\includegraphics[width = .8\linewidth]{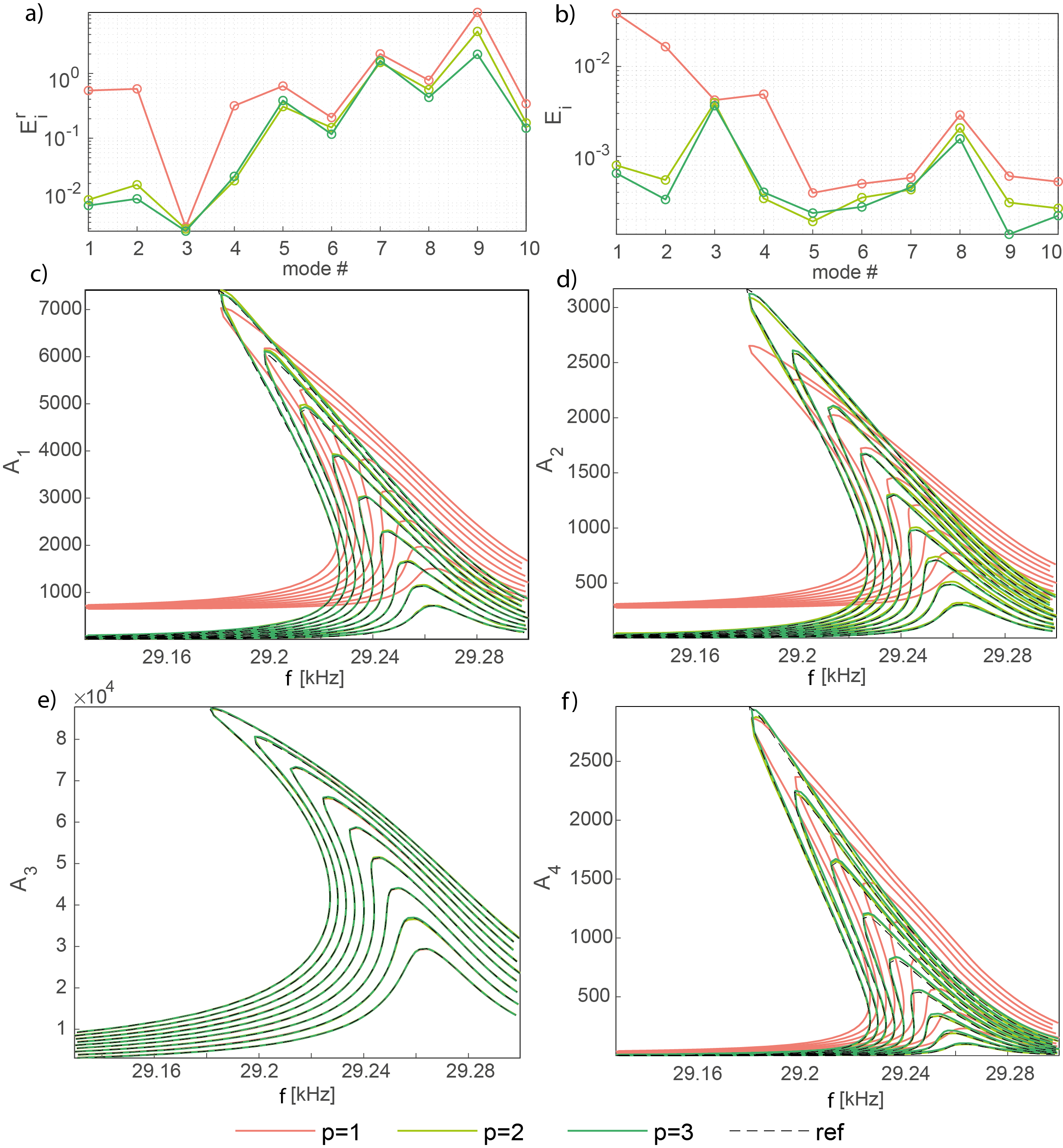}
	\caption{
	Micromirror. Convergence of the analysis increasing the number of reduced variables, with $p=1,2,3$.
	Figures a) and b): error norms $E^r_i$  and $E_i$ defined in Eqs.\eqref{eq:err1}. Figures c)-f): FRFs of the modal amplitudes $A_1-A_4$. The four modes are illustrated in Figure~\ref{fig:perseusmodes} }
	\label{fig:perseusmodal}
\end{figure}

Since the central plate is stiff, we adopt its angle of rotation $\theta$ as reference output  for the FRFs in Figure~\ref{fig:perseusFRF}. It should be noted that the amplitude-induced shift of the peak frequency is small with respect to the absolute value. However, deviations of few Hertzs in working conditions affect the optical performance of projectors; hence, the first requirement for the model is to be highly predictive even of tiny deviations.

This example, that has already been discussed in \cite{DLijnme}, is here revisited using the new arc-length abscissa which corrects the artifacts that were visible at the  edges of the frequency range, and focusing on the convergence of the slave modes and manifolds to the DPIM reference ones.
The DPIM is also used to compute the training data for the DL-ROM method. 
In the following, we consider 63434 snapshots providing 161 samples on each period; the parameter space spanned is  $(\beta,s)$ in $\{1, 1.5, 2, 2.5, 3\}\mu$N $\times [0\!:\!2.0]$. The arc-length abscissa $s$ is rescaled in this example between 0 and 2, and the alignment of the 
data is forced only at peaks of the FRFs, i.e.\ at $s=1$.
One {\it master} mode is used in the DPIM analyses, which implies that the dimension of the reduced space is 2 and the reduced variables can be interpreted as the 
generalised displacement and velocity of the torsional mode. 
The DL-ROM FRFs for the rotation angle are collected in Figure \ref{fig:perseusFRF}
and are very accurate even for $p=1$, i.e,\ with just one reduced variable. Actually, the difference in FRFs 
obtained with increasing $p$ can hardly be appreciated.
The DL-ROM is thus extremely efficient in extracting the key features of the response of this lightly damped device. Indeed, as a first approximation, the mirror can be considered as a conservative system which is well described by a single displacement-like parameter.
\begin{figure}[h]
	\centering
	\includegraphics[width = 0.8\linewidth]{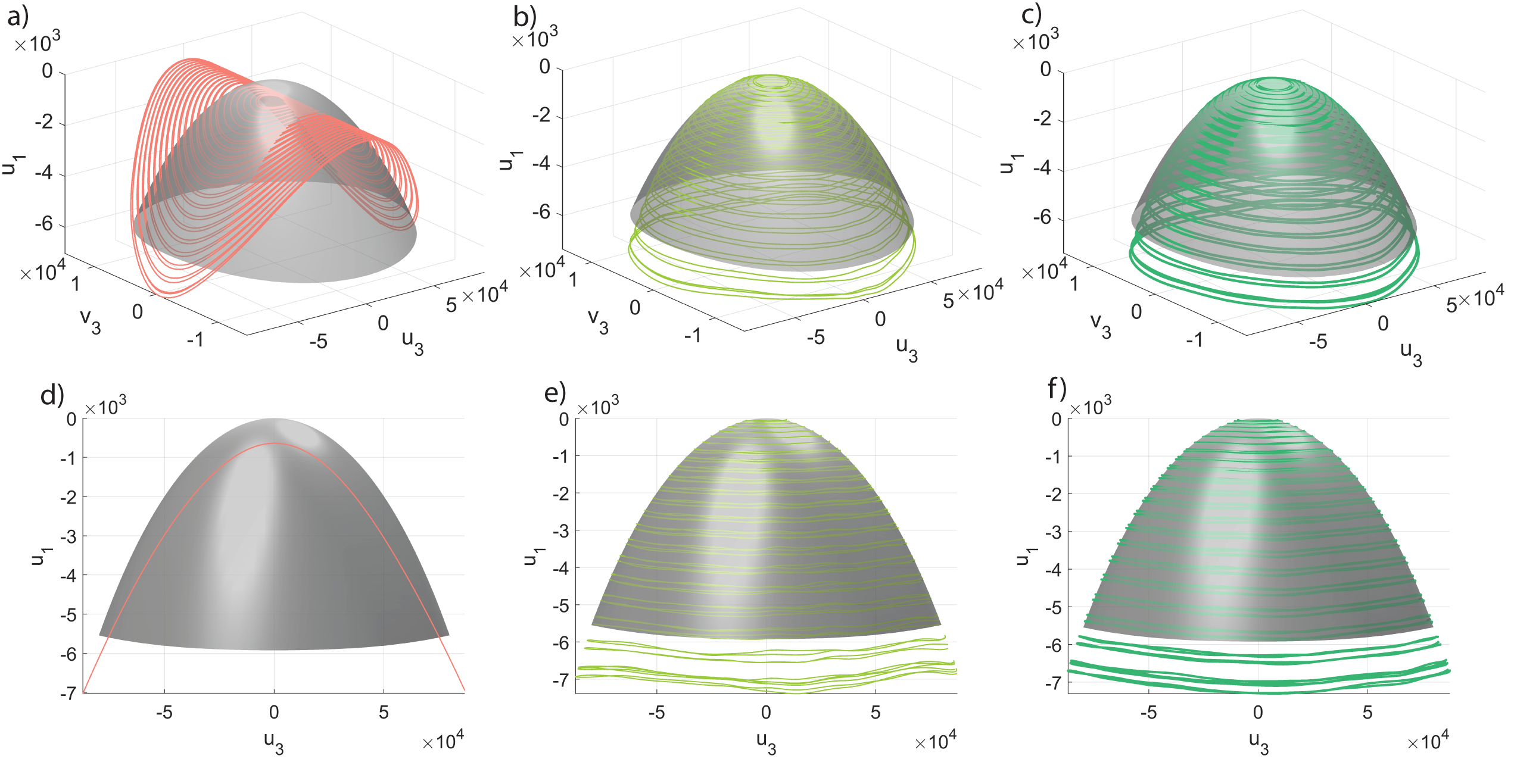}
	\caption{
		Micromirror. Invariant manifolds for the modal displacement $u_1$ in terms of the modal displacement $u_3$ and velocity $v_3$ of the master mode, for $\beta=3.0$. 
		The DPIM manifold is the smooth grey surface, while orbits from the DL-ROM are traced as lines.
		Figures a) and d): different views for $p=1$. Figures b) and e): different views for $p=2$.
		Figures c) and f): different views for $p=3$. It should be noted that the DL-ROM manifolds for 
		$p=1$ are single curvature surfaces.  
	}
	\label{fig:perseusmanifold}
\end{figure} 

The accuracy of the DL-ROM is further investigated in Figure~\ref{fig:perseusmodal} which presents the error measures, Eqs.\eqref{eq:err1}-\eqref{eq:err2}, and the FRFs for the first four modal coordinates $u_i$, Eq.\eqref{eq:umod}, as a function of the forcing frequency. Inspecting the error plots it can be appreciated that, while the master mode is unaffected by the additional reduced variables, these latter on the contrary improve the representation of slave modes. In particular, the low frequency ones benefit from the enrichment of the reduced space, while high-order modes are less sensitive. It should be stressed, anyway, that even if the relative error on these modes remains large, their contribution to the global displacement field is negligible, as illustrated in Figures~\ref{fig:perseusmodal}b).
Moreover, the inclusion of a third reduced variable has a minimal impact on the response, in accordance with the DPIM that perfectly matches the HB FOM with only a single master mode, i.e., two reduced variables \cite{vizzaccaro2021high}. 
Furthermore, the DL-ROM training process is subjected to a certain degree of stochasticity also because of the stochastic gradient descent method employed during the training of neural networks and consequently the results obtained with $p=2$ seem to be, in some conditions, slightly better than the one obtained with $p=3$.
All these remarks are confirmed by the FRFs of the modal coordinates  in Figure \ref{fig:perseusmodal}.  Indeed,  the FRF in Figure \ref{fig:perseusmodal}e), referred to the master mode, is already at convergence with $p=1$, while those in Figure \ref{fig:perseusmodal}c,d,f), plotting the amplitude of the low frequency slave modes, benefit from the introduction of a second reduced variable and are almost insensitive to the third one.
\begin{figure}[h]
	\centering
	\includegraphics[width = 0.8\linewidth]{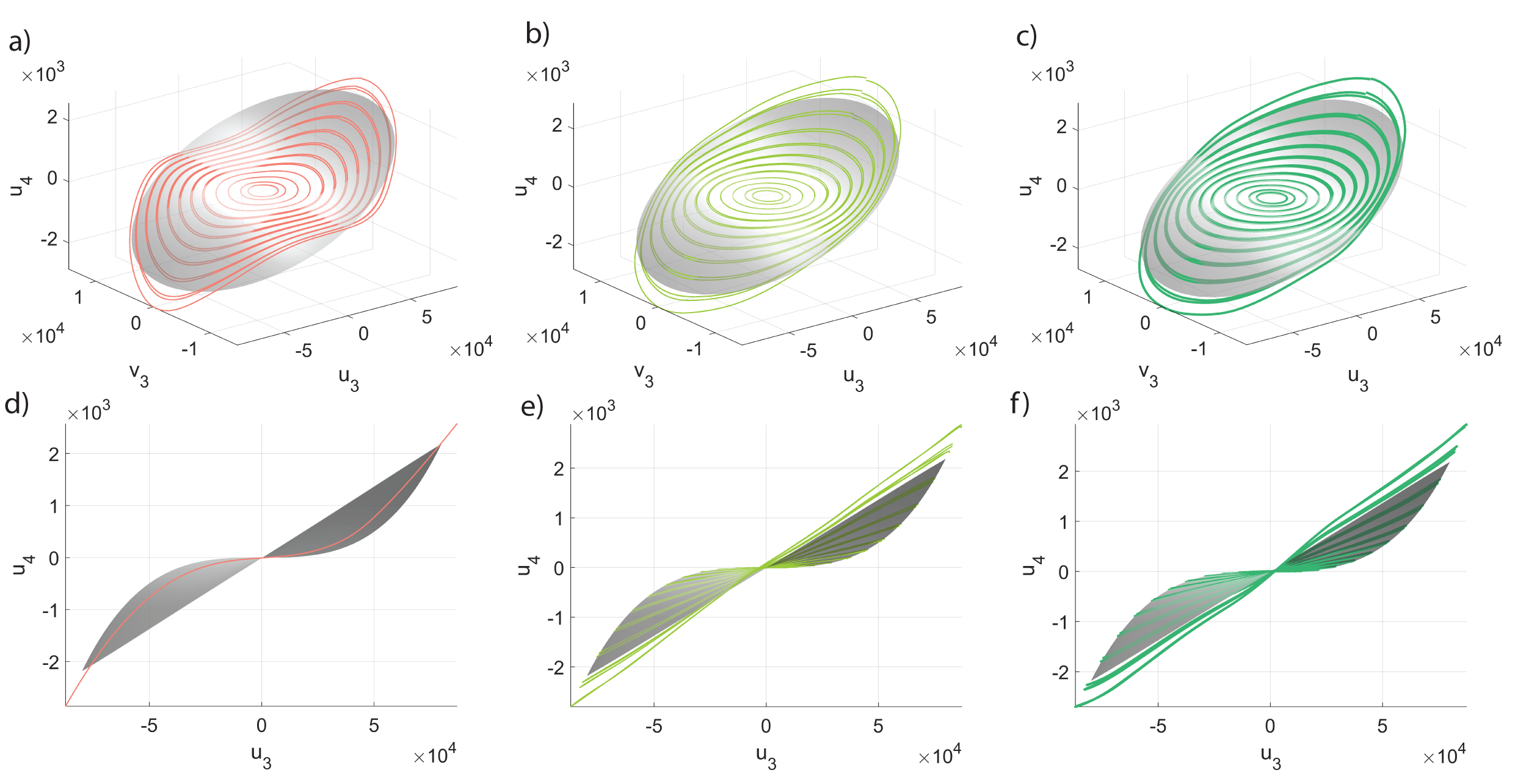}
	\caption{Micromirror. Invariant manifolds for the modal displacement $u_4$ in terms of the modal displacement $u_3$ and velocity $v_3$ of the master mode, for $\beta=3.0$. 
		The DPIM manifold is the smooth grey surface, while orbits from the DL-ROM are traced as lines.
		Figures a) and d): different views for $p=1$. Figures b) and e): different views for $p=2$.
		Figures c) and f): different views for $p=3$. Also in this case the DL-ROM manifolds for 
		$p=1$ are single curvature surfaces.}
	\label{fig:perseusmanifold2}
\end{figure}
In order to further gain insight into the DL-ROM performance, we inspect the invariant manifolds of 
Figure~\ref{fig:perseusmanifold}, computed as explained in the introduction to Section \ref{sec:applications}.
Each manifold plots the modal coordinate $u_i$ 
as a function of $u_3$, the coordinate of the master torsional mode, and of its time derivative $v_3$. 
Two low frequency  slave modes are considered, i.e.\ modes 1 and 4 in Figure~\ref{fig:perseusmodes}. 
The DPIM invariant manifolds are represented as smooth grey surfaces, while DL-ROM orbits are plotted as continuous lines of different colors.
The orbits obtained with $p=1$ are plotted separately in the Figures on the left, in order to better stress the fact that they correspond to single concavity surfaces 
which are, as expected, velocity independent. Indeed, having one single reduced variable, the DL-ROM cannot introduce a dependence on its time derivative.
However, it emerges that
the DL-ROM selects, among all options, the single curvature manifolds that better interpolates the correct one. Manifolds of {\it slave} modes match accurately the DPIM ones already with 2 reduced variables, as clearly  illustrated in Figure~\ref{fig:perseusmodal}. 
\FloatBarrier
\subsection{Internal resonance in a shallow arch}
\label{sec:CCbeam}

Internal resonances (IRs) are associated to energy transfer between modes, and are frequently experienced by MEMS structures mainly due to the very low dissipation. 
Often IRs are strongly linked to the stability of the associated periodic response, and quasi-periodic regimes might arise as a consequence of Neimark-Sacker (NS) bifurcations \cite{mecc21}. The numerical prediction of such phenomena has been tackled only recently
with ad-hoc approaches like the Implicit Condensation \cite{ijnm19} or the DPIM \cite{opreni2022fast}.
Such a task requires a stability analysis which can be run on small ROMs using dedicated continuation tools, and cannot be performed using FOMs, in general.

\begin{figure}[ht]
	\centering
	\includegraphics[width = .8\linewidth]{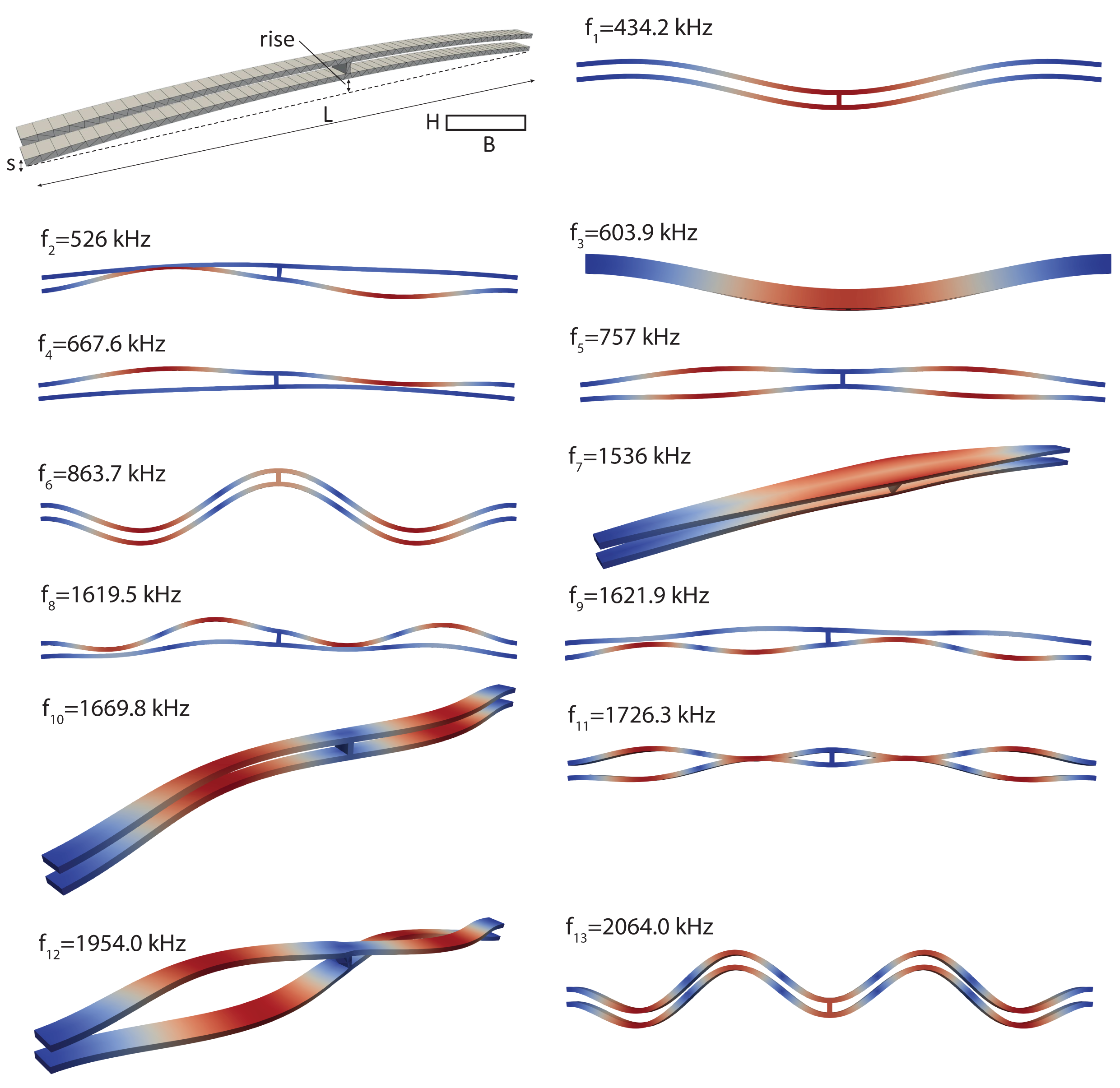}
	\caption{Shallow arch. Geometry and eigenmodes}
	\label{fig:archmodes}
\end{figure}

A typical structure that is prone to 1:2 IR is the shallow arch illustrated in Figure~\ref{fig:archmodes}.
The layout, inspired by the one proposed for a bistable structure in \cite{frangi2015}, has been suitably calibrated so as to trigger the desired IR.
A MEMS resonator presenting the similar geometry and IR has been studied in \cite{zega2022reduced}.
The geometry and the mesh employed for the FOM are illustrated in Figure~\ref{fig:archmodes}. 
The arch dimensions are $B=20\,\mu$m, $H=5\,\mu$m, $L=530\,\mu$m, rise=13.4\,$\mu$m, $s=10\,\mu$m.
The mesh consists of quadratic wedge elements and contains 1971 nodes. The device is made of polycrystalline silicon with density $\rho=2330$\,kg/m$^3$ and
a linear elastic Saint-Venant Kirchhoff constitutive model is assumed, with Young modulus $E=167000$\,MPa and Poisson coefficient $\nu=0.22$ \cite{sharpe1997measurements}. 
The lowest eigenfrequencies of the modelled structure are reported in Figure~\ref{fig:archmodes}.
The quality factor has been set to $Q=500$ and
the actuation is provided by a body force proportional to the first eigenmode $\bfF(t)=\bfM  \bfphi_1 \beta\cos(\omega t)$ with $\beta$ load multiplier.

\begin{figure}[h]
	\centering
	\includegraphics[width = .8\linewidth]{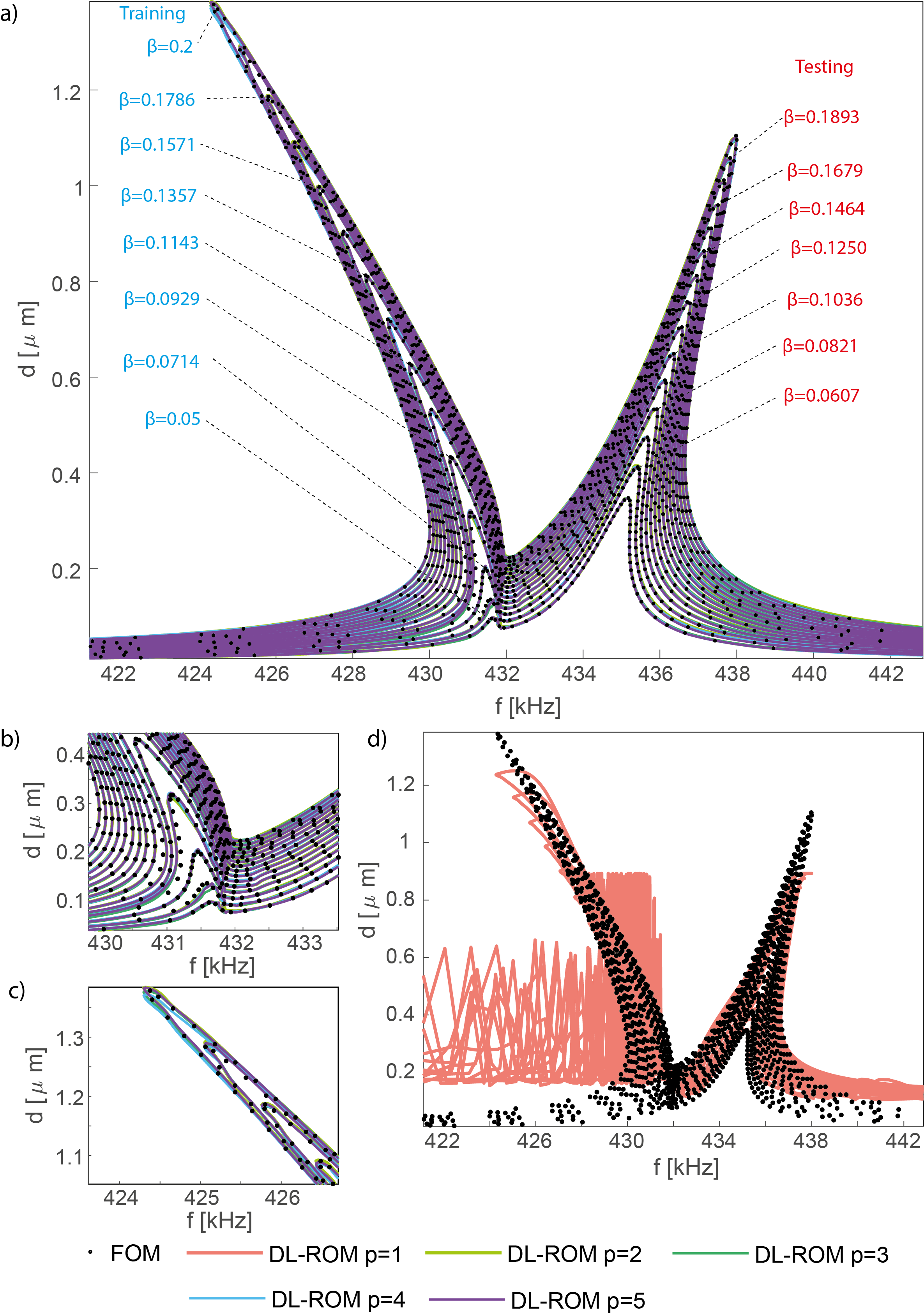}
	\caption{Shallow arch. DLRM FRFs of the mid-span amplitude $d$ and comparison with the reference
		DPIM FRF.
		Figure a): FRFs obtained with 2 and 3 reduced variables. Figures b) and c): details of the FRFs
		in Figure a) showing the excellent accuracy.
		Figure d): FRF obtained with 1 reduced variable. Consistently with the presence of an interaction between to modes, one reduced variable in the DL-ROM cannot describe the full evolution. The FRFs here reported have been obtained with a testing dataset made of 241 339 instances. The inquiry of the DL-ROM is performed in less than 0.5 s using a Tesla V100 32GB GPU and an implementation in the Tensorlow DL framework.}
	\label{fig:arcFRFp23}
\end{figure}

The 1:2 IR occurs as a result of the interaction of the first and the sixth eigenmodes having frequency ratio close to 2 i.e.\ satisfying the necessary condition to induce a 1:2 internal resonance, and leads to qualitative and quantitative changes in the dynamics. 
An in-depth analysis has been developed in Gobat et al. \cite{mecc21}, where 1:2 IR systems are analysed starting their normal  form and the existence of the so-called parabolic modes between the coupled oscillators is demonstrated. Such modes exist on the two branches of the FRF (see Figure \ref{fig:arcFRFp23}) and are associated to the two backbones  that onset from $\omega_1$ and $\omega_6/2$, respectively. 

In order to generate a ROM with the DPIM, the two interacting modes must be both selected as {\it master modes} and the dimension of the reduced space is 4. A critical analysis of the convergence of the DPIM to the full-order HB simulations has been presented in \cite{opreni2022fast}, and on this basis the DPIM results are here considered as an {\it exact} reference.
Addressing this problem with the DL-ROM is a very ambitious task that brings this Deep Learning-based approach to a new level of complexity, as the autoencoder has to automatically recognise the transition between totally different pictures of the dynamical response. Moreover, this example is used to further investigate the
impact of the dimensions of the reduced space for the autoencoder.

\begin{figure}[h]
	\centering
	\includegraphics[width = .8\linewidth]{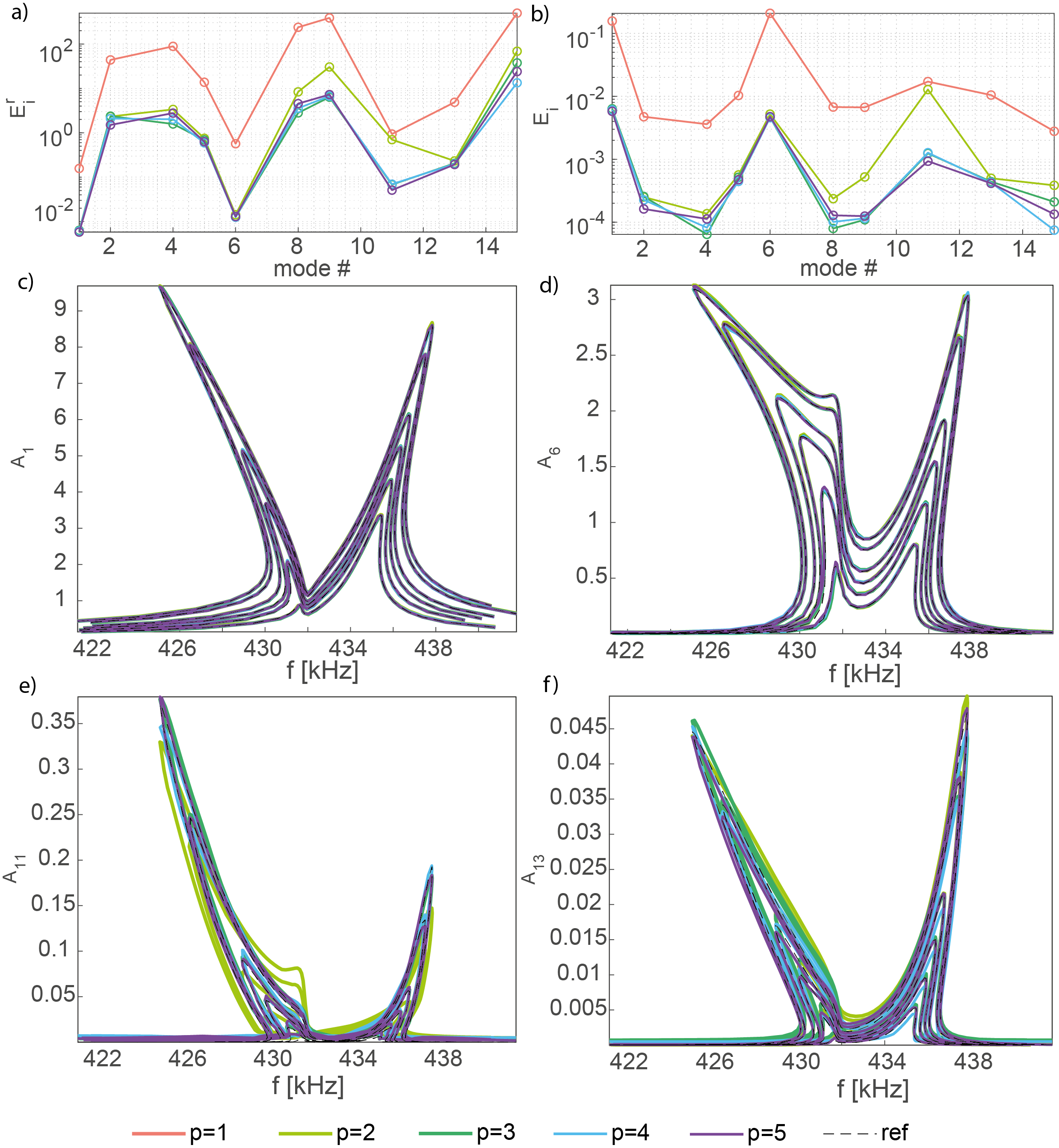}
	\caption{Shallow arch. Convergence of the analysis increasing the number of reduced variables. 
	Figures a)-b): error norms $E^r_i$  and $E_i$, Eqs.\eqref{eq:err1}, 
	for $p=1\shddot 5$. 
	Figures c)-d): FRFs of the master modal amplitudes $A_1,A_6$, $p=2\shddot 5$. 
	Figures e)-f): FRFs of the slave modal amplitudes $A_{11},A_{13}$, $p=2\shddot 5$. 
	For the sake of clarity only the testing $\beta$ values are considered
	and the FRFs for $p=1$ are omitted.}
	\label{fig:FRF_modali_arch}
\end{figure}

Following the same path already considered for the micromirror, in Figure~\ref{fig:arcFRFp23} we first plot the FRFs for the midspan displacement $d$, which is the typical output for this kind of structure.
The reference curves have been computed with the DPIM; note that the DL-ROM FRFs obtained with 2 and 3 reduced variables are almost superposed, thus showing its  convergence, as opposed to the curve obtained with $p=1$.
Indeed, due to the presence of an interaction between two modes, a single reduced variable in the DL-ROM cannot describe the full evolution even for the master modes.

\begin{figure}[h]
	\centering
	\includegraphics[width = .8\linewidth]{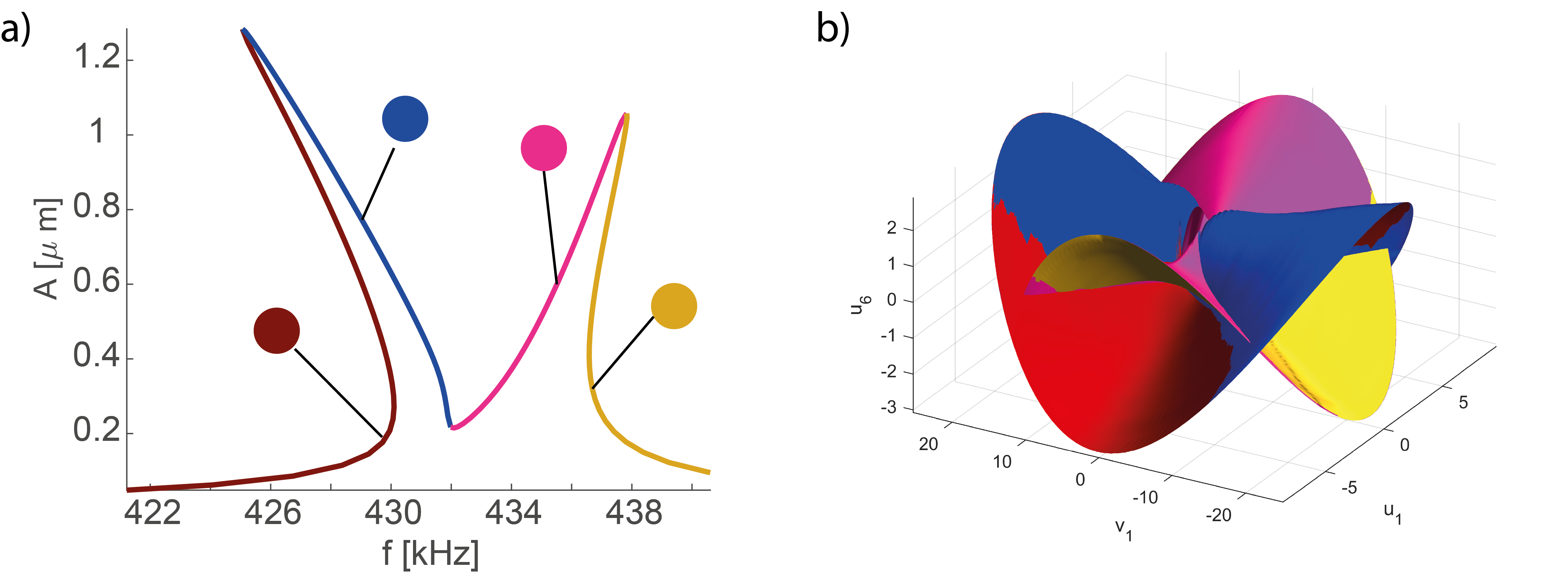}
	\caption{
		Shallow arch. 
		Figure a): reference FRF split in different colour regions. Figure b): reference DPIM envelop of orbits where colours denote the envelops for the regions identified in Figure a). The interaction between the two master modes leads to folding of the surface.}
	\label{fig:refmanifold}
\end{figure}

Next, Figure \ref{fig:FRF_modali_arch} presents the errors,  Eqs.\eqref{eq:err1}-\eqref{eq:err2}, and the FRFs for four selected modal amplitudes, Eq.\eqref{eq:umod}, as a function of the forcing frequency for $p=1\shddot 5$. Inspecting the error plots and the FRFs \ref{fig:FRF_modali_arch}c)-d) for the master modal amplitudes $A_1$ and $A_6$, it can be appreciated that these are already at convergence with $p=2$, consistently with what remarked for the mirror,
where only one {\it master} mode was active and $p=1$ provided accurate FRFs for the master amplitude.

Also in this case, increasing the number of reduced variables has an impact only on the representation of slave modes. The check is performed on modes $11$ and $13$ in Figure~\ref{fig:FRF_modali_arch}, where $p=1\shddot 5$. Convergence to the DPIM develops up to $p=4$, while the inclusion of a fifth reduced variable has a minimal impact  on the response, in accordance with the DPIM that perfectly matches the HB FOM with two modes, i.e.\ four reduced variables \cite{opreni2022fast}.
It should be remarked that these slave amplitudes are at least one order of magnitude smaller than the master ones and, as a consequence, larger errors remain.

\begin{figure}[h]
	\centering
	\includegraphics[width = .8\linewidth]{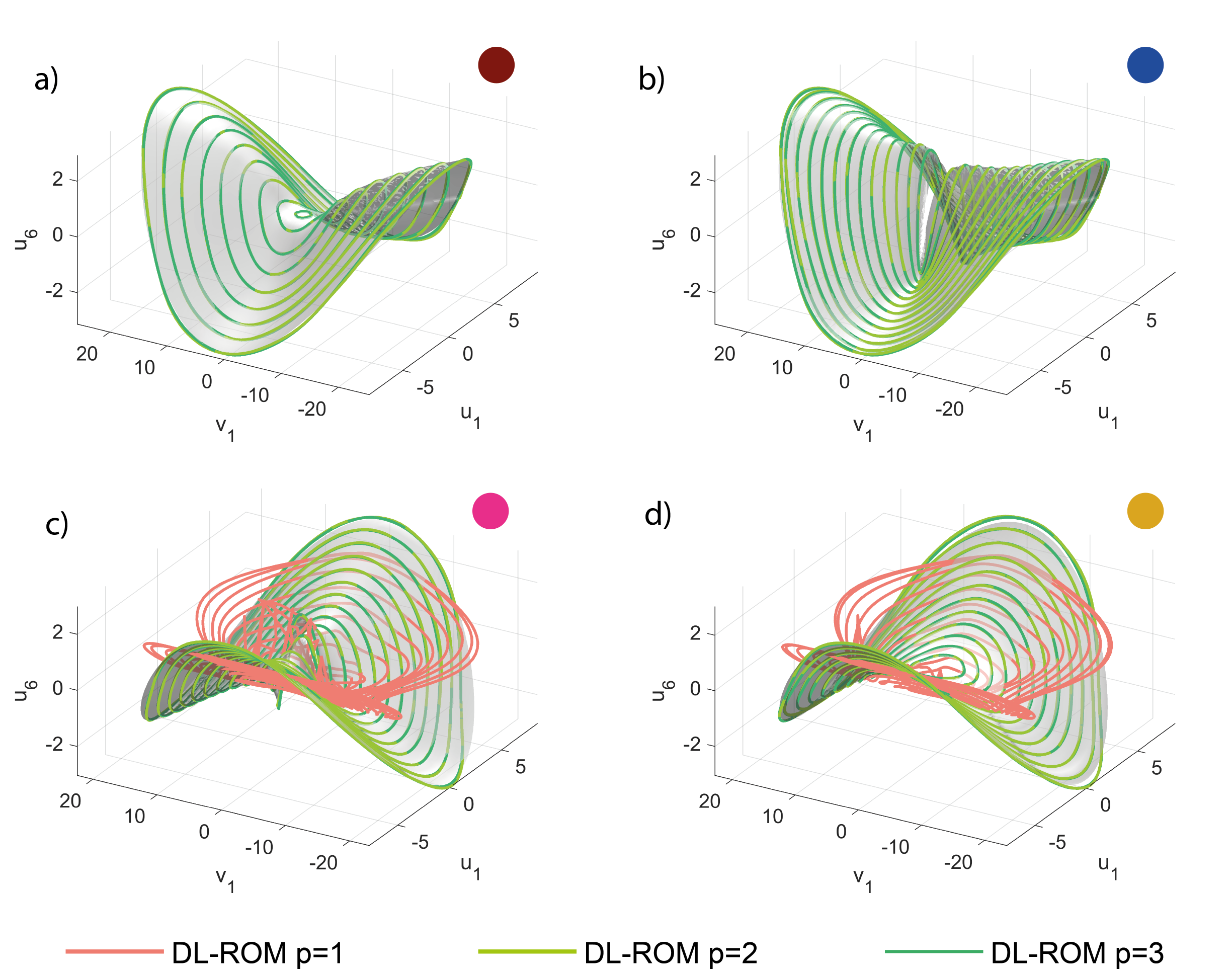}
	\caption{
		Shallow arch. Trajectories $u_6,u_1,v_1$ in a restricted {\it master} space. 
		Figures a)-d): DL-ROM orbits are superposed to the DPIM manifolds for $p=2,3$. The orbits for $p=1$ are only presented on Figures c) and b) and elsewhere omitted for clarity. The 
		{\it master} orbits are already at convergence with $p=2$ (see Fig.~\ref{fig:FRF_modali_arch}).
		The coloured dots refer to the portion of FRF targeted, see Fig. \ref{fig:refmanifold}. }
	\label{fig:arcmanifold}
\end{figure}

In Figure \ref{fig:arcmanifold}, we start comparing the 3D trajectories in the space $u_6,u_1,v_1$, i.e.\ we plot
the modal displacement for the second master mode in terms of the coordinates of the first master.
In this space, the interactions between the two modes leads to foldings and intersections of the orbits envelop as the whole space tends to be filled with orbits.
To tackle these difficulties we will focus on a specific forcing level generating the FRF of Figure \ref{fig:refmanifold}a); furthermore, we will consider five different portions of the FRF denoted by different colors.
The same colors are used to represent the envelops obtained with the DPIM in Figure \ref{fig:refmanifold}b).
In order to better investigate the agreement between DPIM and DL-ROM, the different portions are considered separately in Figures \ref{fig:arcmanifold}a)-d),
where the smooth grey surfaces are the DPIM results, the DL-ROM orbits are traced as lines and the colored dots recall which region of the FRF is being investigated.
From these surfaces, we can appreciate the aforementioned transformation of the manifold which is dominated by a quadratic term turning from negative to positive \cite{mecc21}. Despite the complicated relationship between the orbits on each mode,  the DL-ROM approach provides an almost perfect agreement. The curvatures are well reproduced both qualitatively and quantitatively. 
Since we are considering in Figure \ref{fig:arcmanifold} only master modes, this result confirms our understanding, previously underlined, that the DL-ROM identifies the most critical features of the response and represents them as accurately  as possible: two reduced variables are sufficient to represent their evolution
of two lightly damped, almost conservative modes. 

\begin{figure}[ht]
	\centering
	\includegraphics[width = .8\linewidth]{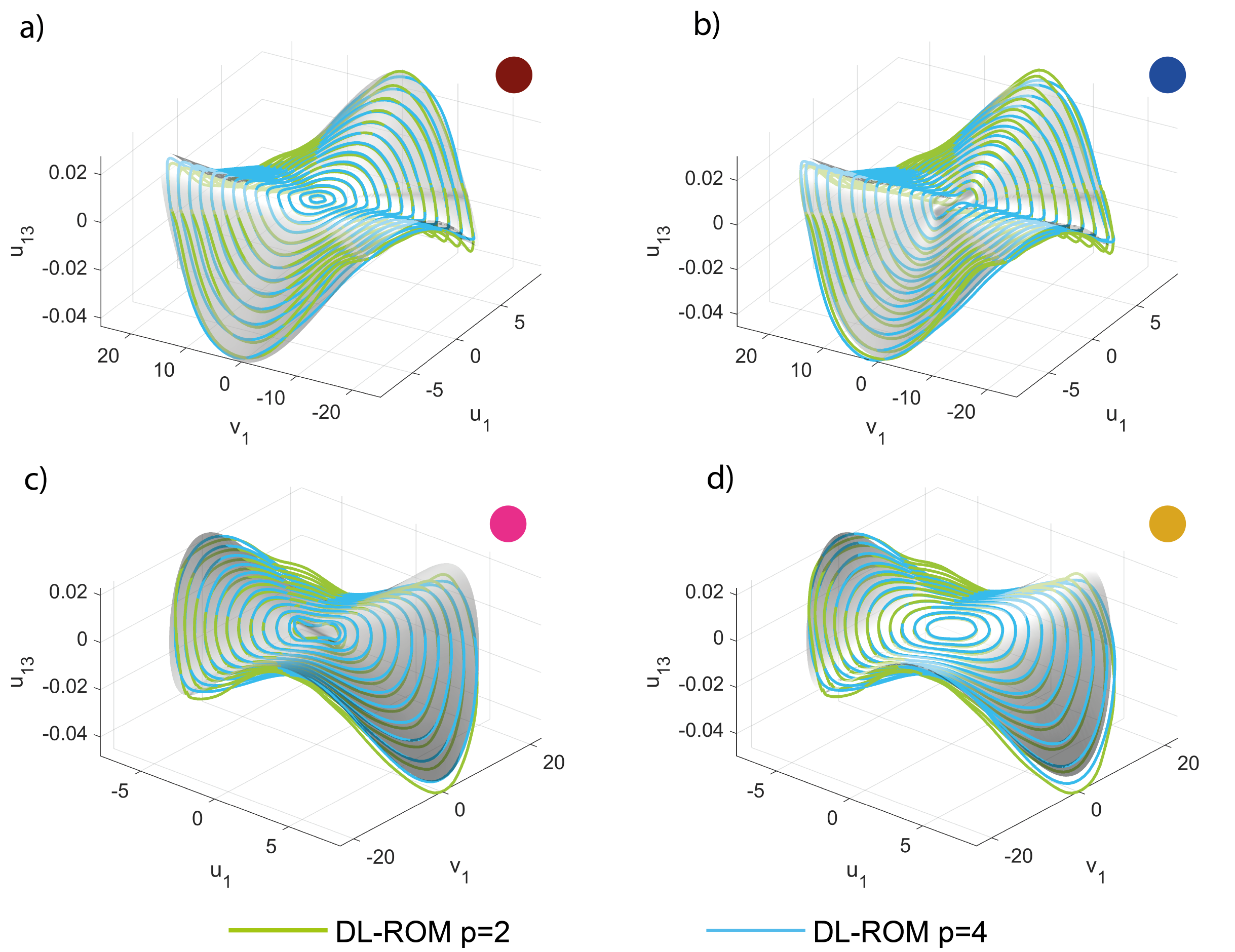}
	\caption{
		Shallow arch. Trajectories $u_{13},u_1,v_1$ in a restricted {\it master} space.
		Figures a)-d): DL-ROM orbits are superposed to the DPIM manifolds for $p=2$ and 4.
		Other $p$ solutions are omitted for sake of clarity and for $p>4$ orbits almost
		coincide (see Fig.~\ref{fig:FRF_modali_arch}). The colored dots refer to the respective FRF region, see Fig. \ref{fig:refmanifold}.}
	\label{fig:arcmanifold2}
\end{figure}
The inspection of the invariant manifolds for slave modes is indeed a much more involved task than for the mirror case, as here the manifolds should be plotted as a function of the four variables $u_1,v_1,u_6,v_6$. 
We thus limit ourselves to tracing envelops of slave trajectories for the modal displacement $u_{13}$ in terms of $u_1,v_1$, illustrated in Figure \ref{fig:arcmanifold2}.
The comparison is done for the same forcing level as before, and the five different regions of the FRF. The smallness of the $u_{13}$ coordinate makes the comparison difficult, as previously 
recalled for the FRFs; however, the convergence of the DL-ROM trajectories to the DPIM smooth surfaces is evident.

\FloatBarrier
\subsection{Electromechanical Disk Resonating Gyroscope}
\label{sec:gyro}

In the previous two examples, an actuation through fictitious body forces has been considered and, thanks to this choice, an extensive validation with
the {\it exact} DPIM approach has been proposed. However, the DPIM has not been applied to multi-physics problems yet, as each extension requires dedicated and costly  developments and coding.

\begin{figure}[h]
	\centering
	\includegraphics[width = .8\linewidth]{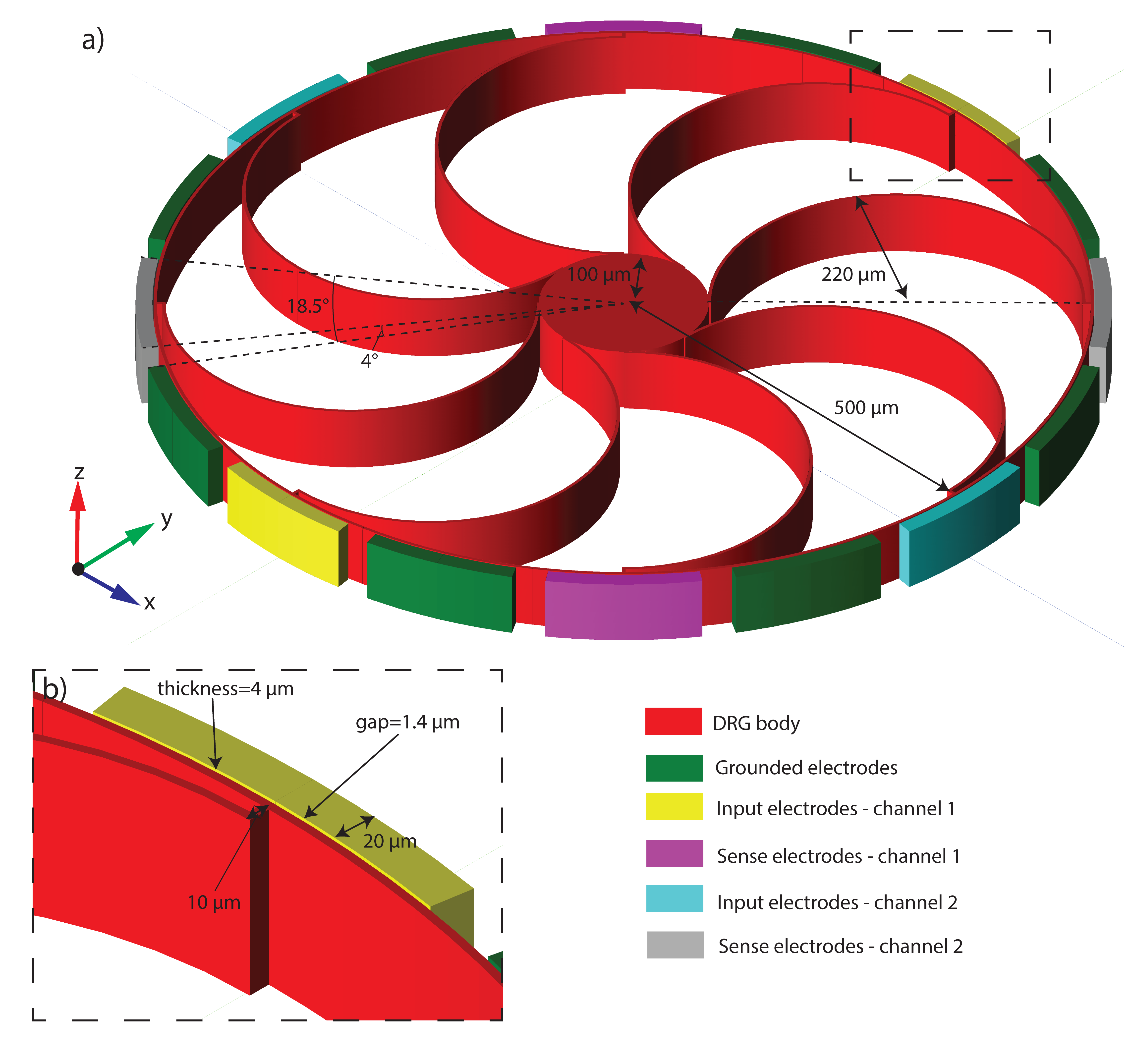}
	\caption{DRG. Figure a): isometric view of the layout. Figure b): details of the radial electrodes}
	\label{fig:geometry_gyro}
\end{figure}
\begin{figure}[h]
\centering
\includegraphics[width = .9\linewidth]{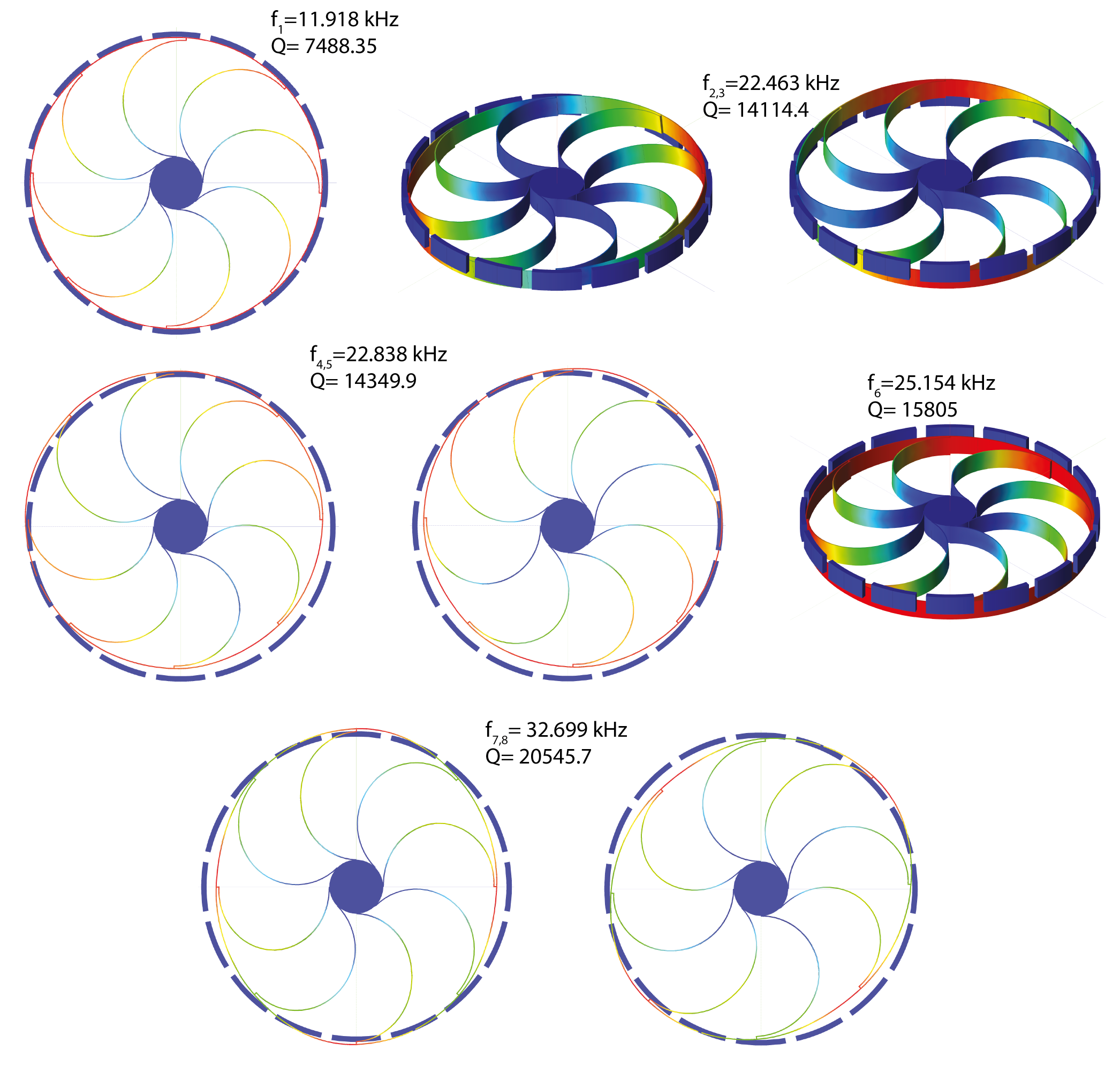}
\caption{DRG: first linear eigenmodes. Due to the symmetry properties, 
many of the eigenmodes come in degenerate couples sharing 
the same eigenfrequency and differing only by a rotation around the $z$ axis.
Modes 7 and 8 are, respectively, the {\it drive} and {\it sense} modes. 
}
\label{fig:modes_gyro}
\end{figure}

Therefore, this application represents a challenging test-bed for the DL-ROM technique, that so far has only been applied to a multi-physics problem in \cite{FrescaManzoniFluids}, where a fluid-structure interaction problem depending on a set of physical parameters has been addressed. We aim at showing the great versatility of the DL-ROM approach, as well as the possibility to easily generalize its construction to face a complex multi-physics problems. Another positive feature, useful in this context, is the non-intrusiveness of DL-ROMs, as snapshots can be readily generated with commercial codes. These capabilities, that are crucial in proposing realistic applications to microstructures, are underlined in this section by addressing an electromechanical Disk Resonanting Gyroscope (DRG) taken from the library of Coventor MEMS+\textsuperscript{TM}\cite{MEMSp}, a leading tool for the analysis of MEMS. The original model available in the software, inspired by the device proposed in \cite{ayazi2001harpss}, has been slightly modified in our benchmarks, and all the details are provided in what follows.

The DRG external ring is modelled with 32 Eulero-Bernulli beams, and each arch suspension is made of 2  Eulero-Bernulli beams. The center support is a rigid body constrained to the ground. 
A series of parallel-plate electrodes are placed around the ring and electrostatic forces are modelled by the software by means of conformal mapping. 
We highlight that, compared to the element size and geometry, the expected displacements are small, so that we neglect geometric nonlinearities and all the nonlinear effects are induced by the electromechanical coupling.

\begin{figure}[h]
	\centering
	\includegraphics[width = .85\linewidth]{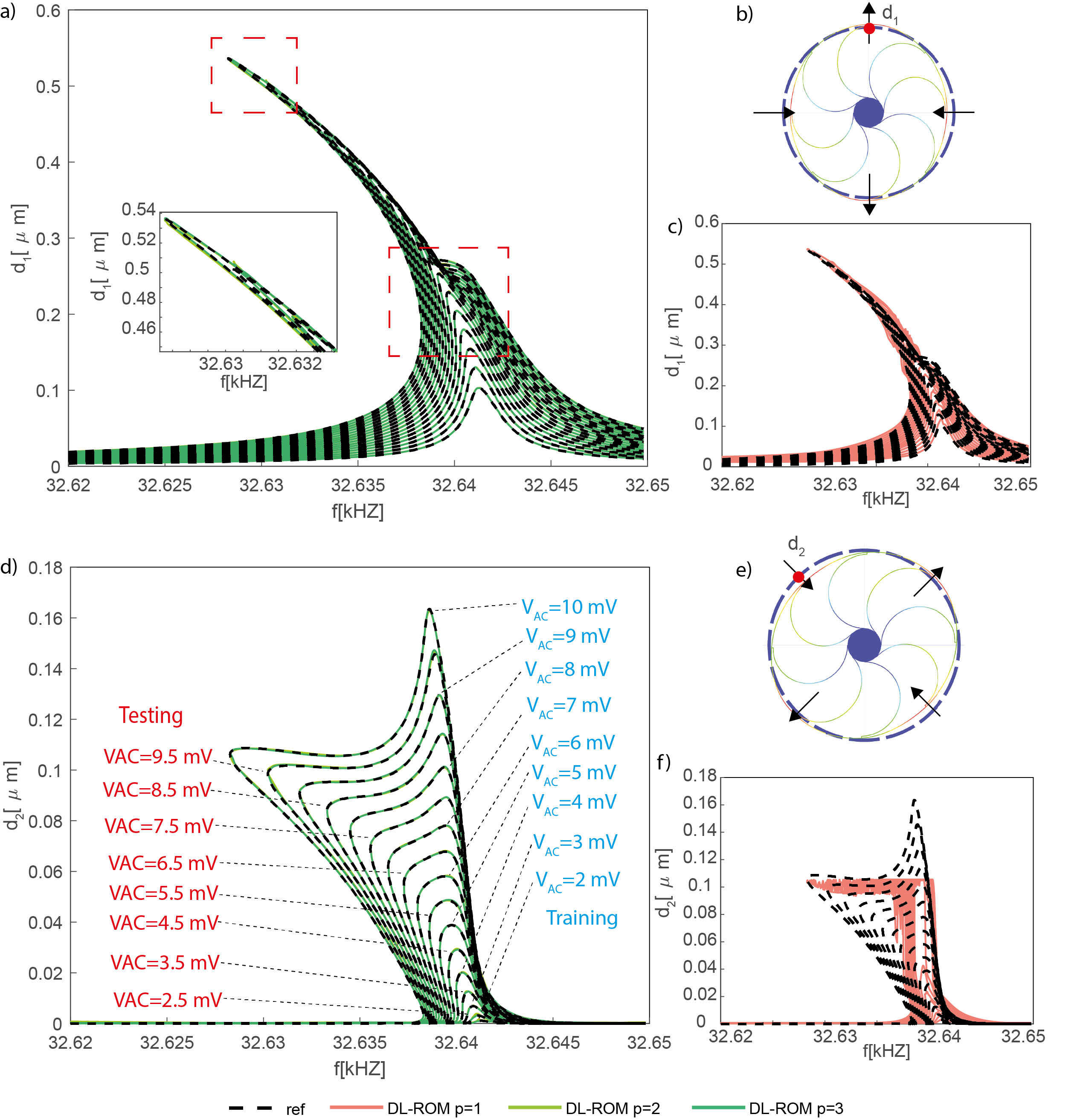}
	\caption{DRG. Figures a) and c): FRFs of the radial displacement of the node indicated by a red circle
	in Figure b), which is representative of the drive mode, i.e.\ mode 7.
	Figure c): comparison between the $p=1$ case and the FOM, clearly highlighting 
	that main dynamic features are not represented adequately with one reduced variable.
	Figure a): comparisons for $p=2,3$.
	The curves below a threshold $V_{\AC}$ correspond to a simple harmonic resonance of the softening drive mode; a plateau starts developing when the sense mode gets autoparametrically activated.
	Data for the sense mode are collected in 
	Figures d) and f) illustrating the FRFs of the radial displacement of the node indicated by a red circle in Figures e), which is representative of the sense mode.
    These FRFs confirm the strong and explosive mode interaction according to which the sense mode reaches abruptly amplitudes of the same order as the drive one. The FRFs here reported have been obtained with a testing dataset made of 1 398 000 instances. The inquiry of the DL-ROM is performed in less than 3 s using a Tesla V100 32GB GPU and an implementation in the Tensorlow DL framework.}
	\label{fig:gyroFRF}
\end{figure}

In Figure \ref{fig:geometry_gyro}, the gyroscope components are coloured according to their function.
This device allows detecting angular velocities around the $z$ axis by exploiting the two degenerate modes 7 and 8 in Figure \ref{fig:modes_gyro} 
which are characterized by radial displacements of the outer ring 
proportional to $\cos(2\theta)$ and $\sin(2\theta)$ with $\theta$ polar coordinate running on the external ring. 
During operations, the drive mode, i.e.\ mode 7, is excited imposing
the bias $V_{\AC}\sin\omega t$ to the blue electrodes 
and the bias $-V_{\AC}\sin\omega t$ to the yellow ones.
The $V_{\AC}$ values range from $2$\,mV up to $10$\,mV guaranteeing an in-plane displacement of approximately 0.54\micr.
A constant potential bias $V_{\DC}=1$\,V is imposed to the gyro 
(red structure). All the remaining electrodes are grounded.
When the device is subjected to an external angular velocity rate $\omega_z$, the Coriolis effect exerts in-plane forces on the ring and activates the sense mode (mode 8 in Figure \ref{fig:modes_gyro}) that is detected by the two couples of independent parallel-plate electrodes
(violet and grey) rotated by 45 degrees with respect to the drive ones. 

However, it has been shown (e.g.\ in \cite{Nitzan15}) that the response of the sense mode largely exceeds predictions based on linear models.
Indeed in the equation of motion for the sense mode, $u_8$,
a term proportional 
to $ u_8 u_7^2V^2_{\DC}\simeq u_8 V^2_{\DC} \sin^2\omega t$ appears
\cite{Polunin17}, 
i.e.\ a linear stiffness term modulated in time at twice the frequency of the mode.
This is responsible for a dynamical instability, called autoparametric resonance \cite{nayfeh95,thomsen2003vibrations}, 
which rather abruptly triggers the activation 
of the sense mode (see Figure \ref{fig:gyroFRF}d).
Simulating the onset of the autoparametric resonance 
is a challenge for any numerical tool
as it induces rapid changes in the overall response, migrating from a single master mode evolution towards a two master modes dynamics. 
For this reason, we assume $\omega_z=0$ in what follows
and focus on autoparametric effects.

\begin{figure}[h]
	\centering
	\includegraphics[width = .95\linewidth]{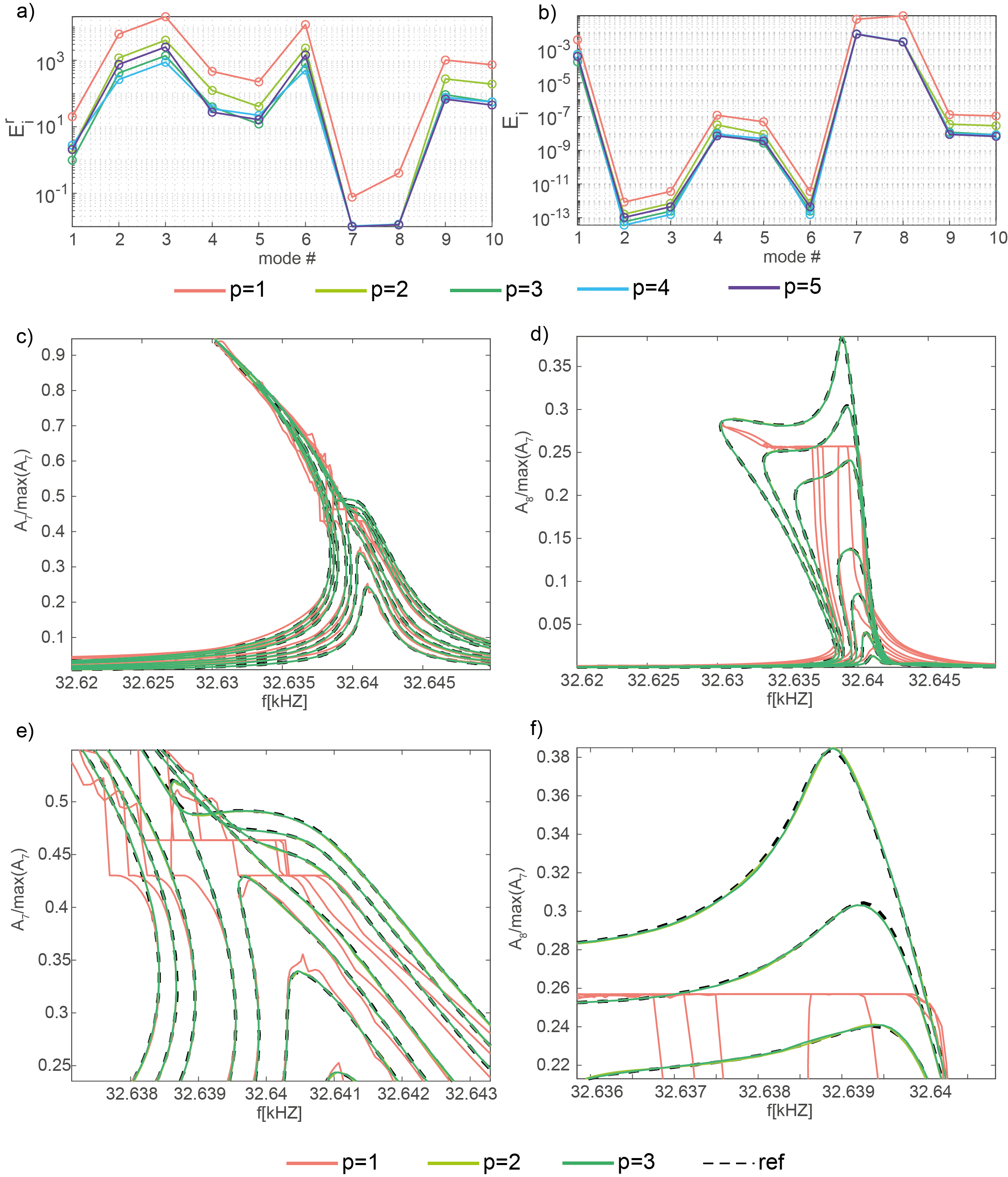}
	\caption{DRG. Convergence of modal coordinates when increasing the number of reduced variables. Figures a) and b): error norms $E^r_i$  and $E_i$, Eqs.\eqref{eq:err1}, for $p=1\shddot 5$. An increase of $p$ over 2 has no effects on the master modes, but improves the representation of slave modes.
	Figures c) and d): FRFs of the master modal amplitudes $A_7,A_8$, $p=1\shddot 3$.
	Figures e) and f): zoom of selected regions in the FRFs}
	\label{fig:FRF_modali_gyro}
\end{figure}
\begin{figure}[h]
	\centering
	\includegraphics[width = .85\linewidth]{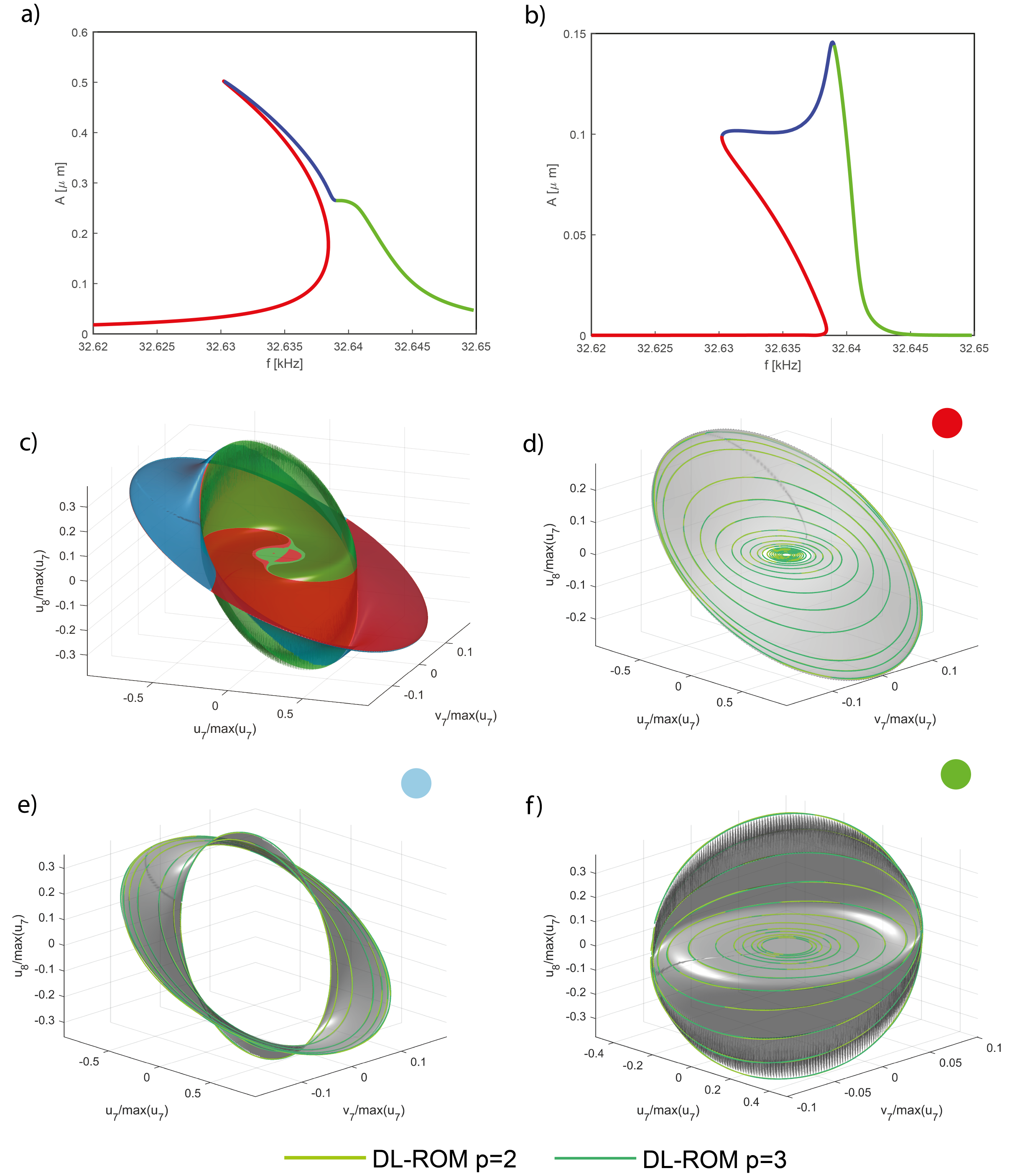}
	\caption{DRG. Envelop of the computed orbits in the 
	$u_8,u_7,v_7$ space for the three branches of the  
	FRFs highlighted in Figure a) and b) (for mode 7 and mode 8, respectively. The whole Conventor MEMS+ envelop is presented in Figure c). 
	The comparison is split in three subregions 
	in Figures d)-f) identified by the colored dots.}
	\label{fig:gyro_manif}
\end{figure}

The DL-ROM has been trained using 149 900 snapshots representing 100 points per period on the parameter space $f=[32.620,32.6498]$\,kHz, $V_{\AC}=[2,3,4,5,6,7,8,9,10]$\,mV. We highlight that the training set is larger than in the previous examples mainly 
due to the changes introduced by the parametric resonance. 
At the onset of autoparametric effects
the drive mode experiences a saturation, while the sense mode abruptly reaches 
comparable amplitudes \cite{avoort10,gallacher2006control}.
The results on a testing set are collected in Figure \ref{fig:gyroFRF} for different choices of number of reduced variables $p$. 
The comparison is performed taking the Coventor results as reference. 

The FRF of the drive amplitude $d_1$, represented by the radial displacement of the node indicated by the red circle in Figures b), is plotted in Figure \ref{fig:gyroFRF}a). 
The curves below a given $V_{\AC}$ denote a simple softening response 
of the drive mode, due to the electrostatic nonlinearities,
while a plateau starts developing when the sense mode gets 
autoparametrically activated.
Inline with the previous examples, in the presence of two master modes
one reduced variable, $p=1$,  cannot represent adequately the main features
and $p=2$ is required. No measurable benefits can be observed when further increasing
$p$. Similar remarks also hold for the FRF of the spurious amplitude $d_2$ plotted in Figure \ref{fig:FRF_modali_gyro}d)-f). 

Also the error estimates \ref{fig:FRF_modali_gyro}a)-b)
show a convergence trend with the increment of $p$ and confirm that 
additional reduced variables only bring benefits to 
slave modes. It is however worth stressing that, in this particular example,
slave modes reach only a negligible amplitudes,  within the numerical errors 
of the procedure, which prevents from presenting meaningful convergence analyses.

Finally, in Figure \ref{fig:gyro_manif} we compare the envelop of orbits in the space $u_8,u_7,v_7$. As for the shallow arch example, 
these master mode orbits tend to fill the whole space. Hence, for the
sake of clarity, we focus on the single FRFs illustrated in 
Figure~\ref{fig:gyro_manif}a),b) which are partitioned in three different regions
generating the Coventor MEMS+ reference envelops of Figure \ref{fig:gyro_manif}c).
The comparison is performed separately in Figures \ref{fig:gyro_manif}d)-f)
for the three different portions highlighted.
Again, $p=1$ is omitted, but the match is excellent already with $p=2$.

\section{Conclusions}
\label{sec:conclusions}

In this work we have focused on two recent computational approaches, namely the Direct Parametrization of Invariant Manifolds (DPIM) and Deep Learning-based ROMs (DL-ROMs), for the efficient numerical simulation of physics-based virtual twins of nonlinear vibrating multi-physics microstructures, providing innovative contributions along two main directions.

On one side, we have brought the applications addressed by POD-DL-ROMs to an unprecedented level by showing that complex dynamical effects can be properly analysed in a very efficient way.
In order to address strongly nonlinear phenomena like internal resonances
and autoparametric effects involving the 
interaction between different modes, we have proposed a new arc-length abscissa which serves as ordering parameter for the collection of snapshots. We have also proved that the POD-DL-ROM construction can be easily extended to multi-physics problems exploiting its non-intrusiveness. This latter indeed allows us to generate snapshots with any FOM, including commercial ones, since FOM solutions are the only simulation data required to train neural network architectures. 

On the other side, for problems with geometric nonlinearities, we have validated the POD-DL-ROM approach against the very recently developed {\it exact} DPIM approach, stressing the striking analogy that exists between the two techniques: namely, both share a nonlinear encoding phase, the generation of a reduced dynamics model, and a likewise nonlinear decoding phase
that allows reconstructing very accurately the full displacement field.
The POD-DL-ROM accurately reproduces the exact invariant manifolds with a minimal set of reduced
variables. 
The problems addressed are characterized by low damping, a distinctive feature of most MEMS working in near-vacuum and these systems can be classified as nearly conservative.
In this context, the POD-DL-ROM approach identifies the major features of the response, i.e., the FRFs together with the trajectories of the $N$ master modes, with only $N$ coordinates.
On the contrary, slave modes might be active in the analysis and their correct reproduction
generally requires to use up to $2N$ reduced variables, as predicted by the DPIM.

The insight gained in this work, associated with the real-time performances of deep learning-based reduced order models, paves the way to the application of the approach in the generation of digital twins of the modelled devices, to be used in the analysis and optimization of the complex systems in which the MEMS are employed, with countless applications in the pervasive Internet of Things. 

\section*{Acknowledgements}
The authors are grateful to Andrea Opreni for the applications of the DPIM and the discussions about invariant manifolds.
GG, SF and AM acknowledge the support of Fondazione Cariplo under grant No. 2019-4608.

\FloatBarrier
\bibliographystyle{unsrt}     

\bibliography{biblio}

\end{document}